\documentclass[10pt,a4paper]{article}
\usepackage{amsthm}
\usepackage{amscd} \usepackage{amsmath}
\usepackage{amsthm}
\usepackage{mathrsfs}
\usepackage{latexsym}
\usepackage{amssymb}
\usepackage{amsfonts}
\theoremstyle{plain}

\title{ { \bf{ Reconstruction of a variety from $\mathscr{O}[[\hbar]]$-modules }}}

\author{\small  By}

\date{\large  Hou-Yi Chen}

\begin{document}

\maketitle

\begin{center}
 $\mathbf{Abstract}$
 \end{center}
We prove that varieties are uniquely determined by the derived
category of $\mathscr{O}[[\hbar]]$-modules with coherent cohomology
which is the same as $\mathscr{O}$-modules proved by A. Bondal and
D. Orlov. We also generalize a theorem of Orlov.\\


\noindent
\emph{{2010} Mathematics subject Classification:}  32C38, 46L65, 53D55.\\

\tableofcontents

\medskip \noindent
\textbf{}

\section{Introduction}
\label{sec.recall}

Let $X$ be a smooth projective variety. In [BO01], A. Bondal and D. Orlov
proved that $X$ can be reconstructed from its bounded derived
category with coherent cohomology when it has ample or anti-ample
canonical bundle. The point of this reconstruction uses the facts that
the Serre functor corresponds to an ample or anti-ample line bundle
(together with a shift).

More precisely, let $X,Y$ be smooth projective varieties, we denote by
$\mathrm{D}^{\mathrm{b}}_{\mathrm{coh}}(\mathscr{O}_{X})$
(resp. $\mathrm{D}^{\mathrm{b}}_{\mathrm{coh}}(\mathscr{O}_{Y})$)
the bounded derived category with coherent cohomology sheaves on
$X$ (resp. $Y$). Then we have the following theorem.\\

\noindent
{\bf{Theorem 1.1.}} ([BO01]). \emph{Let $X$ and $Y$ be two
smooth projective varieties and assume that the (anti)-canonical
bundle of $X$ is ample. If there exists an exact equivalence
$\mathrm{D}^{\mathrm{b}}_{\mathrm{coh}}(\mathscr{O}_{Y})
\stackrel{\sim}\longrightarrow
\mathrm{D}^{\mathrm{b}}_{\mathrm{coh}}(\mathscr{O}_{X})$, then
$X\simeq{Y}$. In particular, the (anti)-canonical bundle of
$Y$ is also ample}.\\

In this paper, we are interested in the study of modules over the ring
$\mathscr{O}_{X}^{\hbar}:=\mathscr{O}[[\hbar]]$ with formal parameter
$\hbar$ and consider the same question as Theorem 1.1. Such modules
naturally appears when studying deformation quantization modules
(DQ-modules) in [KS12]. We obtain the same result for $\mathscr{O}_{X}^{\hbar}$
-modules. More precisely, we prove that\\

\noindent
{\bf{Theorem 1.2.}} \emph{Let $X$ and $Y$ be two
smooth projective varieties and assume that the (anti)-canonical
bundle of $X$ is ample. If there exists an exact equivalence
$\mathrm{D}^{\mathrm{b}}_{\mathrm{coh}}(\mathscr{O}^{\hbar}_{Y})
\stackrel{\sim}\longrightarrow
\mathrm{D}^{\mathrm{b}}_{\mathrm{coh}}(\mathscr{O}_{X}^{\hbar})$, then
$X\simeq{Y}$. In particular, the (anti)-canonical bundle of
$Y$ is also ample}.\\

The method of the proof of Theorem 1.2 is the same as Theorem 1.1 in
[Huy06] and the key point is the existence of the Serre functor.

Recall that for smooth varieties $X,Y$ and $\mathscr{E}\in
\mathrm{D}^{\mathrm{b}}_{\mathrm{coh}}(\mathscr{O}_{X\times{Y}})$,
we have the following Fourier-Mukai transform
\begin{center}
$\Phi_{\mathscr{E}}(-):={\mathrm{R}p_{2}}_{*}(\mathscr{E}\stackrel{
\mathrm{L}}\otimes{p}_{1}^{*}(-)):
\mathrm{D}^{\mathrm{b}}_{\mathrm{coh}}(\mathscr{O}_{X})\rightarrow
\mathrm{D}^{\mathrm{b}}_{\mathrm{coh}}(\mathscr{O}_{Y})$
\end{center}
where $p_{1}:X\times{Y}\rightarrow{X}$
and $p_{2}:X\times{Y}\rightarrow{Y}$  are projections.\\

The following theorem was proved by D. Orlov.\\

\noindent
{\bf{Theorem 1.3.}} ([Orl97]). \emph{Let $X$ and $Y$ be two
smooth projective varieties and let $F:\mathrm{D}^{\mathrm{b}}_{\mathrm{coh}}
(\mathscr{O}_{X})\rightarrow\mathrm{D}^{\mathrm{b}}_{\mathrm{coh}}
(\mathscr{O}_{Y})$ be an equivalence of categories. Then
$F$ is a Fourier-Mukai transform, i.e., there exists
$\mathscr{E}\in\mathrm{D}^{\mathrm{b}}_{\mathrm{coh}}
(\mathscr{O}_{X\times{Y}})$
such that $F\simeq\Phi_{\mathscr{E}}$}.\\

Let $\mathscr{K}\in\mathrm{D}^{\mathrm{b}}_{\mathrm{coh}}
(\mathscr{O}_{X\times{Y}}^{\hbar})$. Define
\begin{center}
(1.1)\hfill$\Phi_{\mathscr{K}}:\mathrm{D}^{\mathrm{b}}_{\mathrm{coh}}
(\mathscr{O}_{X}^{\hbar})\rightarrow\mathrm{D}^{\mathrm{b}}_{\mathrm{coh}}(\mathscr{O}_{Y}^{\hbar})\hfill$
\end{center}
by $\Phi_{\mathscr{K}}(\mathscr{M})={\mathrm{R}p_{2}}_{*}(\mathscr{K}\stackrel{\mathrm{L}}\otimes_{p_{1}^{-1}
\mathscr{O}_{X}^{\hbar}}p_{1}^{-1}\mathscr{M})$,
where $p_{1}:X\times{Y}\rightarrow{X}$
and $p_{2}:X\times{Y}\rightarrow{Y}$  are projections.\\

In this paper, we prove the similar result to Theorem 1.3 as follows:\\

\noindent
{\bf{Theorem 1.4.}} \emph{Let $X$ and $Y$ be two
smooth projective varieties and let $F:\mathrm{D}^{\mathrm{b}}_{\mathrm{coh}}
(\mathscr{O}_{X})\rightarrow\mathrm{D}^{\mathrm{b}}_{\mathrm{coh}}
(\mathscr{O}_{Y})$ be an equivalence of categories. Then
 there exists
$\mathscr{E}^{\mathrm{R}\hbar}\in\mathrm{D}^{\mathrm{b}}_{\mathrm{coh}}
(\mathscr{O}_{X\times{Y}}^{\hbar})$
such that $\Phi_{\mathscr{E}^{\mathrm{R}\hbar}}:
\mathrm{D}^{\mathrm{b}}_{\mathrm{coh}}
(\mathscr{O}_{X}^{\hbar})\rightarrow\mathrm{D}^{\mathrm{b}}_{\mathrm{coh}}(\mathscr{O}_{Y}^{\hbar})$
is an equivalence}.\\

Similarly, we have the well-defined functor as in (1.1) for DQ-modules.
Note that if Question 1.6 below has a positive answer, then Theorem 1.2 may be deduced by the following theorem.\\

\noindent
{\bf{Theorem 1.5.}} ([Pet13, Theorem 3.1.6.]). \emph{Let $X$ and $Y$
be smooth complex projective varieties endowed with
DQ-algebroids $\mathscr{A}_{X}$ and $\mathscr{A}_{Y}$. Let
$\mathscr{K}\in\mathrm{D}^{\mathrm{b}}_{\mathrm{coh}}
(\mathscr{A}_{X\times{Y}})$. Then the following conditions
are equivalent
\begin{itemize}
  \item [(i)] The functor $\Phi_{\mathscr{K}}:
  \mathrm{D}^{\mathrm{b}}_{\mathrm{coh}}(\mathscr{A}_{X})\rightarrow
  \mathrm{D}^{\mathrm{b}}_{\mathrm{coh}}(\mathscr{A}_{Y})$ is fully faithful
  (resp. an equivalence of triangulated categories).
  \item [(ii)] The functor $\Phi_{\mathrm{gr}_{\hbar}\mathscr{K}}:
  \mathrm{D}^{\mathrm{b}}_{\mathrm{coh}}(\mathscr{O}_{X})\rightarrow
  \mathrm{D}^{\mathrm{b}}_{\mathrm{coh}}(\mathscr{O}_{Y})$ is fully faithful
  (resp. an equivalence of triangulated categories).
\end{itemize}}

Therefore, the most interesting question would be the existence of
Fourier-Mukai kernels for DQ-modules.\\

\noindent
{\bf{Question 1.6.}} \emph{Let $\mathscr{A}_{X}$ (resp. $\mathscr{A}_{Y}$) be a DQ-algebroid
on $X$ (resp. $Y$).
If there exists an equivalence of triangulated categories
\begin{center}
$\Phi:\mathrm{D}^{\mathrm{b}}_{\mathrm{coh}}(\mathscr{A}_{X})
\rightarrow\mathrm{D}^{\mathrm{b}}_{\mathrm{coh}}
(\mathscr{A}_{Y})$,
\end{center}
is it a Fourier-Mukai transform}?\\

This paper is organized as follows:
In section 2, we review some results about DQ-modules in [DGS11] and [KS12].
We establish a fully faithful functor between the derived categoy of
algebraic $\mathscr{O}_{X}^{\hbar}$-modules and the derived category of
analytic $\mathscr{O}_{X_{\mathrm{an}}}^{\hbar}$-modules in section 3.
Section 4 contains some results about the triangulated category
$\mathrm{D}^{\mathrm{b}}_{\mathrm{coh}}(\mathscr{O}_{X})$ in [Huy06].
Finally, we prove Theorem 1.2 and Theorem 1.4 in section 5 and
section 6 respectively.\\

\emph{Note:} In this paper, we assume that varieties are smooth over
complex number field $\mathbb{C}$.\\

\begin{center}
$\mathbf{
\mathrm{\mathbf{Notations}}}$
\end{center}
Denote by $\mathbb{K}:=\mathbb{C}$ or $\mathbb{C}^{\hbar}:=\mathbb{C}[[\hbar]]$.\\
Let $(X,\mathscr{A}_{X})$ be a topological space endowed with a sheaf of
algebras over $\mathbb{K}$ or DQ-algebroid $\mathscr{A}_{X}$ over $
\mathbb{C}[[\hbar]]$.\\

We denote by
\begin{itemize}
  \item [(1.1)] $\mathrm{Mod}(\mathscr{A}_{X})$: the abelian category of
  $\mathscr{A}_{X}$-modules,
  \item [(1.2)] $\mathrm{Mod}_{\mathrm{coh}}(\mathscr{A}_{X})$:
  the thick abelian subcategory of $\mathrm{Mod}(\mathscr{A}_{X})$
  consisting of coherent $\mathscr{A}_{X}$-modules,
  \item [(1.3)] $\mathrm{D}(\mathscr{A}_{X})$: the derived category of
  $\mathrm{Mod}(\mathscr{A}_{X})$,
  \item [(1.4)] $\mathrm{D}^{+}(\mathscr{A}_{X})$: the full triangulated
  subcategory consisting of objects with bounded from below cohomology,
  \item [(1.5)] $\mathrm{D}^{\mathrm{b}}(\mathscr{A}_{X})$:
  the full triangulated subcategory consisting of objects with bounded
  cohomology,
  \item [(1.6)] $\mathrm{D}^{\mathrm{b}}_{\mathrm{coh}}
  (\mathscr{A}_{X})$ (resp. $\mathrm{D}^{+}_{\mathrm{coh}}(\mathscr{A}_{X}))$:
  the full triangulated subcategory of $\mathrm{D}^{\mathrm{b}}(\mathscr{A}_{X})$
  (resp. $\mathrm{D}^{+}(\mathscr{A}_{X}))$ consisting of complexes with cohomology
  sheaves belonging to $\mathrm{Mod}_{\mathrm{coh}}(\mathscr{A}_{X})$.
  \item [(1.7)] $(\cdot)^{*}=\mathrm{RHom}_{\mathbb{K}}(\cdot,\mathbb{K})$
  the duality functor on $\mathrm{D}^{\mathrm{b}}(\mathbb{K})$.
\end{itemize}

\noindent
\textbf{Acknowledgments.}
I would like to thank Pierre Schapira
for useful Remarks.

\section{Review on DQ-modules (after [KS12] and [DGS11])}
\label{sec.recall}

In this section, we review some definitions and results in [KS12] and
[DGS11] which we will need in this paper.\\

\noindent
{ $\mathrm{\mathbf{Modules}}$ $\mathrm{\mathbf{over}}$
$\mathrm{\mathbf{formal}}$ $\mathrm{\mathbf{deformations}}$}\\

Let $X$ be a topological space, $\mathscr{A}$ a sheaf of $\mathbb{Z}[\hbar]$-algebra
on $X$ and $\hbar$ a section of $\mathscr{A}$ contained in the center
of $\mathscr{A}$. Let $\mathscr{M}\in\mathrm{Mod}(\mathscr{A})$. Set
\begin{center}
$\widehat{\mathscr{M}}:=\varprojlim\limits_{n}\mathscr{M}/\hbar^{n}\mathscr{M}$.
\end{center}
One says that $\mathscr{M}$ has no $\hbar$-torsion if $\hbar:
\mathscr{M}\rightarrow\mathscr{M}$ is injective and $\mathscr{M}$
is $\hbar$-complete if $\mathscr{M}\rightarrow\widehat{\mathscr{M}}$ is
an isomorphism.

Set
\begin{center}
$\mathscr{A}^{\mathrm{loc}}:=\mathbb{Z}[\hbar,\hbar^{-1}]\otimes_{\mathbb{Z}[\hbar]}
\mathscr{A}$ and $\mathscr{A}_{0}:=\mathscr{A}/\hbar\mathscr{A}$.
\end{center}
We assume the following conditions:
\begin{itemize}
  \item [(i)] $\mathscr{A}$ has no $\hbar$-torsion and is $\hbar$ complete,
  \item [(ii)] $\mathscr{A}_{0}$ is a left Noetherian ring,
  \item [(iii)] there exists a base $\mathfrak{B}$ of open subsets of $X$
  such that for any $U\in\mathfrak{B}$ and any coherent $(\mathscr{A}_{0}|_{U})$-module
  $\mathscr{M}$, we have $H^{n}(U;\mathscr{F})=0$ for any $n>0$.
\end{itemize}

Consider the functor
\begin{center}
$\mathrm{gr}_{\hbar}:\mathrm{D}(\mathscr{A})\rightarrow
\mathrm{D}(\mathscr{A}_{0}),\mathscr{M}\rightarrow\mathrm{gr}_{\hbar}
(\mathscr{M}):=\mathscr{A}_{0}\stackrel{\mathrm{L}}\otimes_{\mathscr{A}}
\mathscr{M}.$
\end{center}
\noindent
{\bf{Lemma 2.1.}} ([KS12, Lemma 1.4.2]). \emph{Let $\mathscr{M}\in\mathrm{D}(\mathscr{A})$ and let
$a\in\mathbb{Z}$. Then we have an exact sequence of $\mathscr{A}_{0}$-modules}
\begin{center}
$0\rightarrow\mathscr{A}_{0}\otimes_{\mathscr{A}}H^{a}(\mathscr{M})\rightarrow
{H}^{a}(\mathrm{gr}_{\hbar}(\mathscr{M}))\rightarrow\mathscr{T}or^{\mathscr{A}}_{1}
(\mathscr{A}_{0},H^{a+1}(\mathscr{M}))\rightarrow{0}$.
\end{center}
${}$\\

For $\mathscr{M},\mathscr{N}\in\mathrm{D}(\mathscr{A})$ and
$\mathscr{P}\in\mathrm{D}(\mathscr{A}^{\mathrm{op}})$, one
has isomorphisms:
\begin{itemize}
  \item [(2.1)] $\mathrm{gr}_{\hbar}(\mathscr{P}\stackrel{\mathrm{L}}\otimes_{\mathscr{A}}\mathscr{M})
  \simeq\mathrm{gr}_{\hbar}(\mathscr{P})\stackrel{\mathrm{L}}\otimes_{\mathscr{A}_{0}}
  \mathrm{gr}_{\hbar}(\mathscr{M})$,
  \item [(2.2)] $\mathrm{gr}_{\hbar}(\mathrm{R}\mathscr{H}om_{\mathscr{A}}
  (\mathscr{M},\mathscr{N}))\simeq\mathrm{R}\mathscr{H}om_{\mathscr{A}_{0}}
  (\mathrm{gr}_{\hbar}\mathscr{M},\mathrm{gr}_{\hbar}\mathscr{N})$.
\end{itemize}
${}$\\
\noindent
{ $\mathrm{\mathbf{Cohomologically}}$ $\hbar$-$\mathrm{\mathbf{complete}}$ $\mathrm{\mathbf{modules}}$}\\

\noindent
{\bf{Definition 2.2.}} One says that an object $\mathscr{M}$ of $\mathrm{D}(\mathscr{A})$
is cohomologically $\hbar$-complete if $\mathrm{R}\mathscr{H}om_{\mathscr{A}}
(\mathscr{A}^{\mathrm{loc}},\mathscr{M})=0$.\\

We denote by $\mathrm{D}_{\mathrm{cc}}(\mathscr{A})$ the full triangulated
subcategory of $\mathrm{D}(\mathscr{A})$ consisting of cohomologically
$\hbar$-complete modules.

Note that $\mathscr{M}\in\mathrm{D}(\mathscr{A})$ is cohomologically $\hbar$-complete
if and only if its image in $\mathrm{D}(\mathbb{Z}[\hbar])$ is cohomologically
$\hbar$-complete.\\

The functor $\mathrm{gr}_{\hbar}$
is conservative on $\mathrm{D}_{\mathrm{cc}}(\mathscr{A})$.\\

\noindent
{\bf{Proposition 2.3.}} ([KS12, Corollary 1.5.9]). \emph{Let $\mathscr{M}\in\mathrm{D}_{\mathrm{cc}}
(\mathscr{A})$. If $\mathrm{gr}_{\hbar}(\mathscr{M})\simeq{0}$, then
$\mathscr{M}\simeq{0}$.}\\

\noindent
{\bf{Corollary 2.4.}} \emph{Let $f:\mathscr{M}\rightarrow\mathscr{N}$ be a morphism
in $\mathrm{D}_{\mathrm{cc}}(\mathscr{A})$. If $\mathrm{gr}_{\hbar}
(f)$ is an isomorphism then $f$ is an isomorphism}.\\

\noindent
{\bf{Theorem 2.5.}} ([KS12, Theorem 1.6.1]). \emph{If $\mathscr{M}\in\mathrm{D}^{\mathrm{b}}_{\mathrm{coh}}
(\mathscr{A})$, then $\mathscr{M}\in\mathrm{D}_{\mathrm{cc}}
(\mathscr{A})$.}\\

The functor $\mathrm{gr}_{\hbar}$ induces the following functor (we keep the same notation):

\begin{center}
(2.3)\hfill$\mathrm{gr}_{\hbar}:\mathrm{D}^{\mathrm{b}}_{\mathrm{coh}}
(\mathscr{A})\rightarrow\mathrm{D}^{\mathrm{b}}_{\mathrm{coh}}
(\mathscr{A}_{0}).\hfill$
\end{center}

\noindent
{\bf{Theorem 2.6.}} ([KS12, Proposition 1.4.5]). \emph{Let $\mathscr{M}\in\mathrm{D}^{\mathrm{b}}_{\mathrm{coh}}
(\mathscr{A})$ and let $a\in\mathbb{Z}$. The conditions below are equivalent:
\begin{itemize}
  \item [(i)] $H^{a}(\mathrm{gr}_{\hbar}(\mathscr{M}))\simeq{0}$,
  \item [(ii)] $H^{a}(\mathscr{M})\simeq{0}$ and $H^{a+1}
  (\mathscr{M})$ has no $\hbar$-torsion.
\end{itemize}}

\noindent
{\bf{Corollary 2.7.}} ([KS12, Corollary 1.4.6]). \emph{The functor $\mathrm{gr}_{\hbar}$ of
(2.3) is conservative.}\\

${}$\\
\noindent
{ $\mathrm{\mathbf{Formal}}$ $\mathrm{\mathbf{extension}}$}\\

Let $\mathscr{R}_{0}$ be a sheaf of rings on a topological space $X$ and let
\begin{center}
$\mathscr{R}:=\mathscr{R}_{0}[[\hbar]]=\prod\limits_{n\geq{0}}
\mathscr{R}_{0}\hbar^{n}$
\end{center}
be the formal extension of $\mathscr{R}_{0}$, whose sections on an open
subset $U$ are formal series $r=\sum_{n=0}^{\infty}r_{j}\hbar^{n}$,
with $r_{j}\in\Gamma(U;\mathscr{R}_{0})$. Consider the associated functor
\begin{center}
$(\cdot)^{\hbar}:\mathrm{Mod}(\mathscr{R}_{0})\rightarrow
\mathrm{Mod}(\mathscr{R})$,
\end{center}
\begin{center}
$\mathscr{N}\mapsto\mathscr{N}[[\hbar]]=\varprojlim\limits_{n}
(\mathscr{R}_{n}\otimes_{\mathscr{R}_{0}}\mathscr{N})$,
\end{center}
where $\mathscr{R}_{n}:=\mathscr{R}/\hbar^{n+1}\mathscr{R}$ is regarded as an
$(\mathscr{R},\mathscr{R}_{0})$-bimodule. Since $\mathscr{R}_{n}$ is free of
finite rank over $\mathscr{R}_{0}$, the functor $(\cdot)^{\hbar}$ is left exact.
We denote by $(\cdot)^{\mathrm{R}\hbar}$ its right derived functor.\\

\noindent
{\bf{Proposition 2.8.}} ([DGS11, Proposition 2.5]). \emph{Let $\mathscr{T}$ be a basis of open subsets of the site
$X$ or, assuming that $X$ is a locally compact topological space, a basis of
compact subsets. Denote by $J_{\mathscr{T}}$ the full subcategory of
$\mathrm{Mod}(\mathscr{R}_{0})$ consisting of $\mathscr{T}$-acyclic objects,
i.e., sheaves $\mathscr{N}$ for which $H^{k}(S;\mathscr{N})=0$ for all
$k>0$ and all $S\in\mathscr{T}$. Then $J_{\mathscr{T}}$ is injective
with respect to the functor $(\cdot)^{\hbar}$. In particular, for
$\mathscr{N}\in{J}_{\mathscr{T}}$, we have $\mathscr{N}^{\hbar}\simeq
\mathscr{N}^{\mathrm{R}\hbar}$.}\\

\noindent
{\bf{Proposition 2.9.}}  \emph{The following sheaves are acyclic for the functor $(\cdot)^{\hbar}$:
\begin{itemize}
  \item [(i)] quasi-coherent modules over the ring $\mathscr{O}_{X}$ of
  regular functions on an algebraic variety $X$,
  \item [(ii)] coherent modules over the ring $\mathscr{O}_{Y}$ of
  holomorphic functions on a complex analytic manifold $Y$.
\end{itemize}}

$\mathrm{\mathbf{Assume}}$ $\mathscr{R}_{0}=\mathscr{A}_{0}$ and $\mathscr{A}=\mathscr{R}_{0}[[\hbar]].$\\

\noindent
{\bf{Proposition 2.10.}} ([DGS11, Proposition 2.8]). \emph{For $\mathscr{N}\in\mathrm{D}^{\mathrm{b}}_{\mathrm{coh}}
(\mathscr{A}_{0})$:
\begin{itemize}
  \item [(i)] there is an isomorphism $\mathscr{N}^{\mathrm{R}\hbar}
\stackrel{\sim}\longrightarrow\mathscr{A}\stackrel{\mathrm{L}}\otimes_{\mathscr{A}_{0}}
\mathscr{N}$,
  \item [(ii)] there is an isomorphism $\mathrm{gr}_{\hbar}(\mathscr{N}^{R\hbar})\simeq
\mathscr{N}$.
\end{itemize}}
In particular, the functor $(\cdot)^{\hbar}$ is exact on $\mathrm{Mod}_{\mathrm{coh}}
(\mathscr{A}_{0})$ and preserves coherence. One thus gets a functor
\begin{center}
(2.4)\hfill$(\cdot)^{\mathrm{R}\hbar}:\mathrm{D}^{\mathrm{b}}_{\mathrm{coh}}
(\mathscr{A}_{0})\rightarrow\mathrm{D}^{\mathrm{b}}_{\mathrm{coh}}(\mathscr{A}).\hfill$
\end{center}

\noindent
{ $\mathrm{\mathbf{DQ}}$-$\mathrm{\mathbf{algebroids}}$}\\

Let $(X,\mathscr{A}_{X})$ be a variety endowed with a DQ-algebroid $\mathscr{A}_{X}$.
Then we have the $\mathbb{C}$-algebroid stack $\mathrm{gr}_{\hbar}(\mathscr{A}_{X})$.

There is a natural functor $\mathscr{A}_{X}\rightarrow\mathrm{gr}_{\hbar}
(\mathscr{A}_{X})$ of $\mathbb{C}$-algebroid stacks. Since $\mathrm{gr}_{\hbar}
\mathscr{A}_{X}\simeq\mathscr{O}_{X}$, this functor induces a functor
\begin{center}
$\iota_{g}:\mathrm{Mod}(\mathscr{O}_{X})\rightarrow\mathrm{Mod}(\mathscr{A}_{X})$.
\end{center}
The functor $\iota_{g}$ admits a left adjoint functor $\mathscr{M}\rightarrow
\mathbb{C}\otimes_{\mathbb{C}^{\hbar}}\mathscr{M}$. The functor $\iota_{g}$ is
exact and fully faithful and it induces a functor
\begin{center}
$\iota_{g}:\mathrm{D}(\mathscr{O}_{X})\rightarrow\mathrm{D}(\mathscr{A}_{X})$.
\end{center}
Moreover, the functor $\mathrm{gr}_{\hbar}$ extends to the algebroid stack $\mathscr{A}_{X}$:
\begin{center}
$\mathrm{gr}_{\hbar}:\mathrm{D}(\mathscr{A}_{X})\rightarrow
\mathrm{D}(\mathscr{O}_{X})$.
\end{center}
Notice that $\mathscr{O}_{X}$ can be endowed with a structure of
$(\mathscr{A}_{X}\otimes\mathscr{A}_{X}^{\mathrm{op}})$-module.
Hence we have
\begin{center}
(2.5)\hfill$\mathrm{gr}_{\hbar}(\mathscr{M})=\mathscr{O}_{X}\stackrel{\mathrm{L}}\otimes
_{\mathscr{A}_{X}}\mathscr{M}$ and $\iota_{g}(\mathscr{M})=_{\mathscr{A}_{X}}\mathscr{O}_{X}
\otimes_{\mathscr{O}_{X}}\mathscr{M}.\hfill$
\end{center}

\noindent
{\bf{Proposition 2.11.}} ([KS12, Proposition 2.3.6]). \emph{The functor $\mathrm{gr}_{\hbar}$ and $\iota_{g}$ define
pairs of adjoint functors $(\mathrm{gr}_{\hbar},\iota_{g})$ and
$(\iota_{g},\mathrm{gr}_{\hbar}[-1])$.}\\

Let $x$ be a closed point in $X$, we denote by $k(x):=\mathbb{C}(x)$ the
skyscraper sheaf supported in $x$ with value $\mathbb{C}$.\\

From the proof of Lemma 2.1, we have the following result.\\

\noindent
{\bf{Proposition 2.12.}} \emph{The functor $\iota_{g}$ induces a
functor
\begin{center}
(2.6)\hfill$\iota_{g}:\mathrm{D}^{\mathrm{b}}_{\mathrm{coh}}(\mathscr{O}_{X})\rightarrow
\mathrm{D}^{\mathrm{b}}_{\mathrm{coh}}(\mathscr{A}_{X}).\hfill$
\end{center}
Moreover, if $\mathscr{F}\in\mathrm{Mod}_{\mathrm{coh}}(\mathscr{O}_{X})$, then
$\mathrm{gr}_{\hbar}\iota_{g}\mathscr{F}\simeq\mathscr{F}[1]\oplus\mathscr{F}$.
In particular, $\mathrm{gr}_{\hbar}\iota_{g}k(x)\simeq{k}(x)[1]\oplus{k}(x)$.}\\

For an algebroid $\mathscr{A}$ on a topological space, we shall write
$``\sigma\in\mathscr{A}"$ instead of $``\sigma\in\mathscr{A}(U)$ for
some open set $U"$.\\

\noindent
{\bf{Definition 2.13.}} An $\mathscr{A}$-module $\mathscr{L}$ is invertible if it
is locally isomorphic to $\mathscr{A}$, namely for any $\sigma\in\mathscr{A}$,
the $\mathscr{E}nd_{\mathscr{A}}(\sigma)$-module $\mathscr{L}(\sigma)$
is locally isomorphic to $\mathscr{E}nd_{\mathscr{A}}(\sigma)$.\\

\noindent
{\bf{Definition 2.14.}} Let $A$ and $A'$ be two sheaves of $\mathbb{C}$-algebras.
An $A\otimes{A}'$-module $L$ is called bi-invertible if there exists locally a
section $\omega$ of $L$ such that $A\ni{a}\mapsto(a\otimes{1})u\in{L}$ and
$A'\ni{a}'\mapsto(1\otimes{a}')u\in{L}$ give  isomorphisms of $A$-modules and
$A'$-modules, respectively.\\

\noindent
{\bf{Definition 2.15.}} For two algebroids $\mathscr{A}$ and
$\mathscr{A'}$, we say that an $(\mathscr{A}\otimes\mathscr{A'})$-module
$\mathscr{L}$ is bi-invertible if for any $\sigma\in\mathscr{A}$ and
$\mathscr{\sigma'}\in\mathscr{A'},\mathscr{L}(\sigma\otimes\sigma')$
is a bi-invertible $\mathscr{E}nd_{\mathscr{A}}(\sigma)\otimes
\mathscr{E}nd_{\mathscr{A'}}(\sigma')$-module.\\

We denote by $X^{a}$ the variety $X$ endowed with the algebroid
$\mathscr{A}_{X}^{\mathrm{op}}$, that is:
\begin{center}
$\mathscr{A}_{X^{a}}=\mathscr{A}_{X}^{\mathrm{op}}$
\end{center}
where the algebroid $\mathscr{A}_{X}^{\mathrm{op}}$ is defined by
$\mathscr{A}_{X}^{\mathrm{op}}(U)=(\mathscr{A}_{X}(U))^{\mathrm{op}}$
($U$ is open in $X$).  Note that this is also a DQ-algebroid.\\

\noindent
{\bf{Proposition 2.16.}} ([KS12, Lemma 2.4.1]). \emph{Let $\mathscr{M}$ be an $(\mathscr{A}_{X}\otimes
\mathscr{A}_{X^{a}})$-module. Then the following conditions are equivalent.}

\begin{itemize}
   \item [(i)] \emph{$\mathscr{M}$ is a bi-invertible $(\mathscr{A}_{X}\otimes
\mathscr{A}_{X^{a}})$-module (see definition 2.15),}
  \item [(ii)] \emph{$\mathscr{M}$ is invertible as an $\mathscr{A}_{X}$-module, that is,
$\mathscr{M}$ is locally isomorphic to $\mathscr{A}_{X}$ as an $\mathscr{A}_{X}$-module
(see definition 2.13),}
  \item [(iii)] \emph{$\mathscr{M}$ is invertible as an $\mathscr{A}_{X^{a}}$-module}.
\end{itemize}

\noindent
{\bf{Proposition 2.17.}} ([KS12]). \emph{The category of bi-invertible $(\mathscr{A}_{X}\otimes
\mathscr{A}_{X^{a}})$-modules has a structure of a tensor category with unit
object $\mathscr{A}_{X}$.}\\

\noindent
{\bf{Corollary 2.18.}} \emph{For a bi-invertible $(\mathscr{A}_{X}\otimes
\mathscr{A}_{X^{a}})$-module $\mathscr{P}$, the functor
\begin{center}
(2.7)\hfill$\mathscr{P}\stackrel{\mathrm{L}}\otimes_{\mathscr{A}_{X}}(\cdot):\mathrm{D}^{\mathrm{b}}_{\mathrm{coh}}
(\mathscr{A}_{X})\longrightarrow\mathrm{D}^{\mathrm{b}}_{\mathrm{coh}}
(\mathscr{A}_{X})\hfill$
\end{center}
is an equivalence with inverse functor
\begin{center}
(2.8)\hfill$\mathscr{P}^{\otimes{-1}}\stackrel{\mathrm{L}}\otimes_{\mathscr{A}_{X}}(\cdot):
\mathrm{D}^{\mathrm{b}}_{\mathrm{coh}}
(\mathscr{A}_{X})\longrightarrow\mathrm{D}^{\mathrm{b}}_{\mathrm{coh}}
(\mathscr{A}_{X})\hfill$
\end{center}
where $\mathscr{P}^{\otimes{-1}}:=\mathrm{R}\mathscr{H}om_{\mathscr{A}_{X}}
(\mathscr{P},\mathscr{A}_{X})$.}\\

\noindent
{\bf{Corollary 2.19.}} \emph{For a bi-invertible $(\mathscr{A}_{X}\otimes\mathscr{A}_{X^{a}})$-
module $\mathscr{P}$, $\mathrm{gr}_{\hbar}\mathscr{P}$ is a line bundle.}\\

\noindent
{\bf{Remark 2.20.}} The above results hold in particular when replacing $\mathscr{A}_{X}$
by $\mathscr{O}_{X}[[\hbar]]$.\\

\section{A fully faithful functor}
\label{sec.recall}

Let $X$ be a smooth algebraic variety with associated complex manifold
$X_{\mathrm{an}}$. We have a natural morphism
\begin{center}
(3.1)\hfill$\iota:X_{\mathrm{an}}\rightarrow{X}\hfill$
\end{center}
of topological spaces and a morphism $\iota^{-1}\mathscr{O}_{X}\rightarrow
\mathscr{O}_{X_{\mathrm{an}}}$ of sheaves of rings. Hence we get a functor
\begin{center}
$(\cdot)_{\mathrm{an}}:\mathrm{Mod}_{\mathrm{coh}}(\mathscr{O}_{X})\longrightarrow
\mathrm{Mod}_{\mathrm{coh}}(\mathscr{O}_{X_{\mathrm{an}}})$
\end{center}
which induces the functor of full triangulated subcategories formed by coherent
cohomology (we keep the same notation)
\begin{center}
(3.2)\hfill$(\cdot)_{\mathrm{an}}:\mathrm{D}^{\mathrm{b}}_{\mathrm{coh}}
(\mathscr{O}_{X})\longrightarrow\mathrm{D}^{\mathrm{b}}_{\mathrm{coh}}
(\mathscr{O}_{X_{\mathrm{an}}}).\hfill$
\end{center}
Set
\begin{center}
$\mathscr{O}_{X}^{\hbar}:=\mathscr{O}_{X}[[\hbar]]$ and
$\mathscr{O}_{X_{\mathrm{an}}}^{\hbar}:=\mathscr{O}_{X_{\mathrm{an}}}[[\hbar]]$.
\end{center}
We also have the following functors induced by (2.3) and (2.4):
\begin{center}
(3.3)\hfill$\mathrm{gr}_{\hbar}:\mathrm{D}^{\mathrm{b}}_{\mathrm{coh}}
(\mathscr{O}_{X}^{\hbar})$ (resp. $\mathrm{D}^{\mathrm{b}}_{\mathrm{coh}}
(\mathscr{O}_{X_{\mathrm{an}}}^{\hbar}))\longrightarrow
\mathrm{D}^{\mathrm{b}}_{\mathrm{coh}}(\mathscr{O}_{X})$
(resp. $\mathrm{D}^{\mathrm{b}}_{\mathrm{coh}}
(\mathscr{O}_{X_{\mathrm{an}}}))\hfill$
\end{center}
and
\begin{center}
(3.4)\hfill$(\cdot)^{\mathrm{R}\hbar}:\mathrm{D}^{\mathrm{b}}_{\mathrm{coh}}
(\mathscr{O}_{X_{\mathrm{an}}})\longrightarrow\mathrm{D}^{\mathrm{b}}_{\mathrm{coh}}
(\mathscr{O}_{X_{\mathrm{an}}}^{\hbar}).\hfill$
\end{center}
\noindent
{\bf{Definition 3.1.}} We define the functor $(\cdot)^{\hbar}_{\mathrm{an}}$ to be
the compositions of the functors (3.2), (3.3) and (3.4):
\begin{center}
(3.5)\hfill $(\cdot)^{\hbar}_{\mathrm{an}}:\mathrm{D}^{\mathrm{b}}_{\mathrm{coh}}
(\mathscr{O}_{X}^{\hbar})\stackrel{\mathrm{gr}_{\hbar}}\longrightarrow
\mathrm{D}^{\mathrm{b}}_{\mathrm{coh}}(\mathscr{O}_{X})\stackrel{(\cdot)_{\mathrm{an}}}
\longrightarrow\mathrm{D}^{\mathrm{b}}_{\mathrm{coh}}(\mathscr{O}_{X_{\mathrm{an}}})
\stackrel{(\cdot)^{\mathrm{R}\hbar}}\longrightarrow\mathrm{D}^{\mathrm{b}}_{\mathrm{coh}}
(\mathscr{O}_{X_{\mathrm{an}}}^{\hbar}).\hfill$
\end{center}

\noindent
{\bf{Theorem 3.2.}} \emph{Assume $X$ is projective. Then the functor $(\cdot)^{\hbar}_{\mathrm{an}}$
of (3.5) is fully faithful.}\\

\noindent
{\emph{{Proof}.}} For $\mathscr{M},\mathscr{N}\in\mathrm{D}^{\mathrm{b}}_{\mathrm{coh}}
(\mathscr{O}_{X}^{\hbar})$, we have the following natural morphism
\begin{center}
(3.6)\hfill$\mathrm{RHom}_{\mathscr{O}_{X}^{\hbar}}(\mathscr{M},\mathscr{N})\longrightarrow
\mathrm{RHom}_{\mathscr{O}_{X_{\mathrm{an}}}^{\hbar}}((\mathscr{M})^{\hbar}_{\mathrm{an}},
(\mathscr{N})^{\hbar}_{\mathrm{an}}).\hfill$
\end{center}
Taking $\mathrm{gr}_{\hbar}$ to (3.6) and using Serre's GAGA theorem (cf. [Che10] Corollary 1.4),
(2.2) and Proposition 2.10, we get the following isomorphism
\begin{center}
$\mathrm{RHom}_{\mathscr{O}_{X}}(\mathrm{gr}_{\hbar}\mathscr{M},
\mathrm{gr}_{\hbar}\mathscr{N})\stackrel{\sim}\longrightarrow
\mathrm{RHom}_{\mathscr{O}_{X_{\mathrm{an}}}}((\mathrm{gr}_{\hbar}\mathscr{M})_{\mathrm{an}},
(\mathrm{gr}_{\hbar}\mathscr{N})_{\mathrm{an}})$.
\end{center}
Since the functor $\mathrm{gr}_{\hbar}$ is conservative by Corollary 2.7, this implies
that (3.6) is an isomorphism and the result follows.$\hfill \square$\\

We denote by $\Omega_{X}^{p}$ (resp. $\Omega_{X_{\mathrm{an}}}^{p})$) the
sheaf of regular (resp. holomorphic) $p$-forms on $X$ (resp. $X_{\mathrm{an}}$)
and set $\Omega_{X}:=\Omega_{X}^{\mathrm{dim}X}$ (resp. $\Omega_{X_{\mathrm{an}}}:=
\Omega_{X_{\mathrm{an}}}^{\mathrm{dim}_{X_{\mathrm{an}}}}$) and set
$\Omega_{X}^{\hbar}:=\Omega_{X}[[\hbar]]$ and $\Omega_{X_{\mathrm{an}}}^{\hbar}:
=\Omega_{X_{\mathrm{an}}}[[\hbar]]$. Then we have\\

\noindent
{\bf{Lemma 3.3.}} \emph{We have $(\Omega_{X}^{\hbar})^{\hbar}_{\mathrm{an}}=
\Omega^{\hbar}_{X_{\mathrm{an}}}$}.

\section{Review on the triangulated category $\mathrm{D}^{\mathrm{b}}_{\mathrm{coh}}
(\mathscr{O}_{X})$}
\label{sec.recall}
In this section, we review some results in [Huy06].\\

\noindent
{\bf{Definition 4.1.}} A collection $\Xi$ of objects in a triangulated
category $\mathcal{D}$ is a spanning class of $\mathcal{D}$ if for all
$B\in\mathcal{D}$ the following two conditions hold:
\begin{itemize}
  \item [(i)] If $\mathrm{Hom}(A,B[i])=0$ for all $A\in\Xi$ and
all $i\in\mathbb{Z}$ then $B\simeq{0}$.
  \item [(ii)] If $\mathrm{Hom}(B[i],A)=0$ for all $A\in\Xi$ and all
$i\in\mathbb{Z}$, then $B\simeq{0}$.
\end{itemize}

\noindent
{\bf{Proposition 4.2.}} ([Huy06, Proposition 3.17]). \emph{Let $X$ be a projective variety. Then the objects of
the form $k(x)$ with $x\in{X}$ a closed point span the triangulated category
$\mathrm{D}^{\mathrm{b}}_{\mathrm{coh}}(\mathscr{O}_{X})$.}\\

\noindent
{\bf{Definition 4.3.}} Let $\mathcal{D}$ be a $\mathbb{C}$-linear
triangulated category with a Serre functor $S$. An object $P\in\mathcal{D}$
is called point like of codimension $d$ if
\begin{itemize}
  \item [(i)] $S(P)\simeq{P}[d]$
  \item [(ii)] $\mathrm{Hom}(P,P[i])=0$ for $i<0$ and
  \item [(iii)] $k(P):=\mathrm{Hom}(P,P)\simeq\mathbb{C}$.
\end{itemize}
An object $P$ satisfying (iii) is called simple.\\

\noindent
{\bf{Proposition 4.4.}} ([Huy06, Lemma 4.5]). \emph{Let $X$ be a projective variety. Suppose $\mathscr{F}$
is a simple object (i.e., $\mathrm{Hom}(\mathscr{F},\mathscr{F})\simeq\mathbb{C})$
in $\mathrm{D}^{\mathrm{b}}_{\mathrm{coh}}(\mathscr{O}_{X})$ with zero
dimension support. If $\mathrm{Hom}(\mathscr{F},\mathscr{F}[i])=0$ for $i<0$, then
\begin{center}
$\mathscr{F}\simeq{k}(x)[m]$
\end{center}
for some closed point $x\in{X}$ and some integer $m$.}\\

\noindent
{\bf{Proposition 4.5.}} ([Huy06, Proposition 4.6]). \emph{Let $X$ be a projective variety. Suppose that $\Omega_{X}$
or $\Omega_{X}^{*}$ is ample. Then the point like objects in
$\mathrm{D}^{\mathrm{b}}_{\mathrm{coh}}(\mathscr{O}_{X})$ are the objects which
are isomorphic to $k(x)[m]$, where $x\in{X}$ is a closed point and $m\in\mathbb{Z}$.}\\

\noindent
{\bf{Definition 4.6.}} Let $\mathcal{D}$ be a triangulated category with a Serre
functor $S$. An object $L\in\mathcal{D}$ is called invertible if for any point like
object $P\in\mathcal{D}$ there exists $n_{P}\in\mathbb{Z}$ (depending also on $L$)
such that
\[ \mathrm{Hom}(L,P[i]) = \left\{
\begin{array}{ll}
k(P) & \mbox{if } i=n_{P} \\[2pt]
0 & \mbox{otherwise.}
\end{array} \right. \]

\noindent
{\bf{Proposition 4.7.}} ([Huy06, Proposition 4.9]). \emph{Let $X$ be a projective variety. Any invertible object in
$\mathrm{D}^{\mathrm{b}}_{\mathrm{coh}}(\mathscr{O}_{X})$ is of the form
$L[m]$ with $L$ a line bundle on $X$ and $m\in\mathbb{Z}$. Conversely,
if $\Omega_{X}$ or $\Omega_{X}^{*}$ is ample, then for any line bundle
$L$ and any $m\in\mathbb{Z}$ the object
$L[m]\in\mathrm{D}^{\mathrm{b}}_{\mathrm{coh}}(\mathscr{O}_{X})$
is invertible.}\\

\noindent
{\bf{Lemma 4.8.}} ([Huy06, Lemma 7.2]). \emph{Let $X$ be a projective variety, $x\in{X}$ a closed point, and
$\mathscr{F}\in\mathrm{D}^{\mathrm{b}}_{\mathrm{coh}}(\mathscr{O}_{X})$.
Suppose $\mathrm{Hom}(\mathscr{F},k(y)[i])=0$ for any closed point
$y\neq{x}$ and any $i\in\mathbb{Z}$ and $\mathrm{Hom}(\mathscr{F},
k(x)[i])=0$ for $i<0$ or $i>\mathrm{dim}X$. Then
$\mathscr{F}$ is isomorphic to a sheaf concentrated in $x\in{X}$.}\\

\section{Proof of Theorem 1.2}
\label{sec.recall}

In this section, we shall prove Theorem 1.2 as stated in the introduction.\\

Let $(X,\mathscr{A}_{X})$ be a projective variety endowed with a DQ-algebroid
$\mathscr{A}_{X}$. Recall from [Che10] that the morphism $\iota:X_{\mathrm{an}}
\rightarrow{X}$ of (3.1) induces an analytic DQ-algebroid $\mathscr{A}_{X_{\mathrm{an}}}$
on $X_{\mathrm{an}}$ and the natural functor
\begin{center}
$(\cdot)_{\mathrm{an}}:\mathrm{Mod}(\mathscr{A}_{X})\rightarrow
\mathrm{Mod}(\mathscr{A}_{X_{\mathrm{an}}})$
\end{center}
induces the following equivalence
\begin{center}
(5.1)\hfill$(\cdot)_{\mathrm{an}}:\mathrm{D}^{\mathrm{b}}_{\mathrm{coh}}
(\mathscr{A}_{X})\stackrel{\sim}\longrightarrow\mathrm{D}^{\mathrm{b}}_{\mathrm{coh}}
(\mathscr{A}_{X_{\mathrm{an}}}).\hfill$
\end{center}
First we show that the existence of the Serre functor on
$\mathrm{D}^{\mathrm{b}}_{\mathrm{coh}}(\mathscr{A}_{X_{\mathrm{an}}})$
induces a Serre functor on $\mathrm{D}^{\mathrm{b}}_{\mathrm{coh}}
(\mathscr{A}_{X})$ by (5.1).\\

\noindent
{\bf{Theorem 5.1.}} \emph{Let $(X,\mathscr{A}_{X})$ be a
projective variety endowed with a DQ-algebroid $\mathscr{A}_{X}$.
Then the triangulated category $\mathrm{D}^{\mathrm{b}}_{\mathrm{coh}}
(\mathscr{A}_{X})$ has a Serre functor defined by
$S_{X}:\mathscr{M}\rightarrow\Omega_{X}^{\mathscr{A}}\otimes_{\mathscr{A}_{X}}
\mathscr{M}[\mathrm{dim}X]$ where $\Omega_{X}^{\mathscr{A}}$ is an
invertible $\mathscr{A}_{X}$-module such that
$(\Omega_{X}^{\mathscr{A}})_{\mathrm{an}}\simeq\Omega_{X_{\mathrm{an}}}^{\mathscr{A}}$
where $\Omega_{X_{\mathrm{an}}}^{\mathscr{A}}$ is defined in [KS12] such that
$\mathrm{gr}_{\hbar}\Omega_{X_{\mathrm{an}}}^{\mathscr{A}}\simeq
\Omega_{X_{\mathrm{an}}}$.}\\

\noindent
{\emph{{Proof}.}} By the equivalence of (5.1), there exists a coherent
$\mathscr{A}_{X}$-module $\Omega_{X}^{\mathscr{A}}$ such that
$(\Omega_{X}^{\mathscr{A}})_{\mathrm{an}}\simeq\Omega_{X_{\mathrm{an}}}^{\mathscr{A}}$.
By [KS12, Corollary 3.3.4], we have for $\mathscr{M},\mathscr{N}\in
\mathrm{D}^{\mathrm{b}}_{\mathrm{coh}}(\mathscr{A}_{X})$
\begin{center}
$\mathrm{RHom}_{\mathscr{A}_{X_{\mathrm{an}}}}
(\mathscr{N}_{an},\Omega_{X_{\mathrm{an}}}^{\mathscr{A}}
\stackrel{\mathrm{L}}\otimes_{\mathscr{A}_{X_{\mathrm{an}}}}\mathscr{M}_{\mathrm{an}})
[\mathrm{dim}(X_{\mathrm{an}})]\stackrel{\sim}\longrightarrow
(\mathrm{RHom}_{\mathscr{A}_{X_{\mathrm{an}}}}(\mathscr{M}_{\mathrm{an}},
\mathscr{N}_{\mathrm{an}}))^{*}$.
\end{center}
Since
\begin{center}
$\mathrm{RHom}_{\mathscr{A}_{X}}(\mathscr{N},\Omega_{X}^{\mathscr{A}}
\stackrel{\mathrm{L}}\otimes_{\mathscr{A}_{X}}\mathscr{M})[\mathrm{dim}(X)]
\simeq\mathrm{RHom}_{\mathscr{A}_{X_{\mathrm{an}}}}
(\mathscr{N}_{\mathrm{an}},\Omega_{X_{\mathrm{an}}}^{\mathscr{A}}
\stackrel{\mathrm{L}}\otimes_{\mathscr{A}_{X_{\mathrm{an}}}}\mathscr{M}_{X_{\mathrm{an}}})
[\mathrm{dim}(X_{\mathrm{an}})]$
\end{center}
and
\begin{center}
$(\mathrm{RHom}_{\mathscr{A}_{X}}(\mathscr{M},\mathscr{N}))^{*}\simeq
(\mathrm{RHom}_{\mathscr{A}_{X_{\mathrm{an}}}}(\mathscr{M}_{\mathrm{an}},
\mathscr{N}_{\mathrm{an}}))^{*}$,
\end{center}
the result follows. Since $\mathrm{gr}_{\hbar}(\Omega_{X}^{\mathscr{A}})_{\mathrm{an}}\simeq
(\mathrm{gr}_{\hbar}\Omega_{X}^{\mathscr{A}})_{\mathrm{an}}\simeq\Omega_{X_{\mathrm{an}}}$,
$\mathrm{gr}_{\hbar}\Omega_{X}^{\mathscr{A}}\simeq\Omega_{X}$ a locally free
$\mathscr{O}_{X}$-module of rank 1 and hence $\Omega_{X}^{\mathscr{A}}$ has no
$\hbar$-torsion by Theorem 2.6. Hence $\Omega_{X}^{\mathscr{A}}$ is invertible.$\hfill \square$\\

\noindent
{\bf{Corollary 5.2.}} \emph{Let $X$ be a projective variety. Then the triangulated
category $\mathrm{D}^{\mathrm{b}}_{\mathrm{coh}}(\mathscr{O}_{X}^{\hbar})$
has a Serre functor defined by $S_{X}:\mathscr{M}\rightarrow
\Omega_{X}^{\hbar}\otimes_{\mathscr{O}_{X}^{\hbar}}\mathscr{M}
[\mathrm{dim}(X)]$ where $\Omega_{X}^{\hbar}$ is an invertible
$\mathscr{O}_{X}^{\hbar}$-module such that $(\Omega_{X}^{\hbar})^{\hbar}_{\mathrm{an}}
\simeq\Omega_{X_{\mathrm{an}}}^{\hbar}$ (cf. Lemma 3.3).}\\

\noindent
{\emph{{Proof}.}} The proof is the same as Theorem 5.1 and using Theorem 3.2.$\hfill \square$\\

We recall the following definition.\\

\noindent
{\bf{Definition 5.3.}} We say that two functors $F,F':\mathcal{A}\rightarrow
\mathcal{B}$ from category $\mathcal{A}$ to $\mathcal{B}$ are isomorphic if
there exists a morphism of functors $\varphi:F\rightarrow{F}'$ such that for
any object $A\in\mathcal{A}$ the induced morphism $\varphi_{A}:
F(A)\rightarrow{F}'(A)$ is an isomorphism in $\mathcal{B}$.\\

\noindent
{\bf{Proposition 5.4.}} \emph{Let $X$ and $Y$ be projective varieties. Assume that
$F:\mathrm{D}^{\mathrm{b}}_{\mathrm{coh}}(\mathscr{A}_{X})\stackrel{\sim}\longrightarrow
\mathrm{D}^{\mathrm{b}}_{\mathrm{coh}}(\mathscr{A}_{Y})$ is an equivalence.
Let $S_{X}$ (resp. $S_{Y}$) be the Serre functor of $X$ (resp. $Y$) in
Theorem 5.1. Then there exists an isomorphism}
\begin{center}
$F\circ{S}_{X}\simeq{S}_{Y}\circ{F}$.
\end{center}
\noindent
{\emph{{Proof}.}} We denote by $\mathcal{A}:=\mathrm{D}^{\mathrm{b}}_{\mathrm{coh}}(\mathscr{A}_{X})$
and $\mathcal{B}:=\mathrm{D}^{\mathrm{b}}_{\mathrm{coh}}(\mathscr{A}_{Y})$.
Since $F$ is fully faithful, one has for any two
objects $A,B\in\mathcal{A}$
\begin{center}
$\mathrm{RHom}(A,S_{X}(B))\simeq\mathrm{RHom}(F(A),F(S_{X}B))$
\end{center}
and
\begin{center}
$\mathrm{RHom}(B,A)\simeq\mathrm{RHom}(F(B),F(A))$.
\end{center}
Together with the two isomorphisms
\begin{center}
$\mathrm{RHom}(A,S_{X}(B))\simeq\mathrm{RHom}(B,A)^{*}$
\end{center}
and
\begin{center}
$\mathrm{RHom}(F(B),F(A))\simeq\mathrm{RHom}(F(A),S_{Y}F(B))^{*}$,
\end{center}
this yields a functorial isomorphism
\begin{center}
$\mathrm{RHom}(F(A),FS_{X}(B))\simeq
\mathrm{RHom}(F(A),S_{Y}F(B))$.
\end{center}
Using the hypothesis that $F$ is an equivalence and in particular that any
object in $\mathcal{B}$ is isomorphic to some $F(A)$ one concludes that
there exists a functor isomorphism $F\circ{S}_{X}\simeq{S}_{Y}\circ{F}$.$\hfill \square$\\

\noindent
{\bf{Lemma 5.5.}} \emph{Let $X$ be a projective variety. Then the objects of the
form $\iota_{g}k(x)$ with $x\in{X}$ a closed point span the
triangulated category $\mathrm{D}^{\mathrm{b}}_{\mathrm{coh}}(\mathscr{A}_{X})$.}\\

\noindent
{\emph{{Proof}.}} First we prove the condition (i) in Definition 4.1:
$\mathrm{Hom}_{\mathscr{A}_{X}}(\mathscr{M},\iota_{g}k(x))[i]=0$, then we get
\begin{center}
$\mathrm{Hom}_{\mathscr{A}_{X}}(\mathscr{M},\iota_{g}k(x))[i]=
\mathrm{Hom}_{\mathscr{O}_{X}}(\mathrm{gr}_{\hbar}\mathscr{M},
k(x))[i]=0$ for all $x\in{X}$.
\end{center}
Since objects $k(x)$ span $\mathrm{D}^{\mathrm{b}}_{\mathrm{coh}}
(\mathscr{O}_{X})$ by Proposition 4.2, we get $\mathrm{gr}_{\hbar}\mathscr{M}
\simeq{0}$ which implies $\mathscr{M}\simeq{0}$ by Corollary 2.7. The
condition (ii) is similar to (i) by Proposition 2.11.$\hfill \square$\\

\noindent
{\bf{Proposition 5.6.}} \emph{We have $\mathrm{Hom}_{\mathscr{A}_{X}}
(\iota_{g}k(x),\iota_{g}k(x))\simeq\mathrm{Hom}_{\mathscr{O}_{X}}
(\mathrm{gr}_{\hbar}\iota_{g}k(x),k(x))\simeq
\mathrm{Hom}_{\mathscr{O}_{X}}(k(x),\mathrm{gr}_{\hbar}\iota_{g}
k(x)[-1])\simeq\mathrm{Hom}_{\mathscr{O}_{X}}(k(x),k(x))\simeq\mathbb{C}$.}\\

\noindent
{\emph{{Proof}.}} Consider the following distinguished triangle
\begin{center}
$k(x)[1]\longrightarrow{k}(x)[1]\oplus{k}(x)\longrightarrow
{k}(x)\stackrel{0}\longrightarrow{k}(x)[2]$
\end{center}
which induces the exact sequence
\begin{center}
$0\longrightarrow\mathrm{Hom}(k(x),k(x))\longrightarrow
\mathrm{Hom}(k(x)[1]\oplus{k}(x),k(x))\longrightarrow{0}$
\end{center}
since $\mathrm{Hom}(k(x)[i],k(x))=0$ for $i>0$. This implies
that $\mathrm{Hom}(k(x)[1]\oplus{k}(x),k(x))\simeq
\mathrm{Hom}(k(x),k(x))\simeq\mathbb{C}$.

Now the results follow from Proposition 2.11 and Proposition 2.12.$\hfill \square$\\

\noindent
{\bf{Lemma 5.7.}} There exists a natural isomorphism
\begin{center}
$S_{X}(\iota_{g}k(x))\simeq\iota_{g}k(x)[\mathrm{dim}X]$.
\end{center}
\noindent
{\emph{{Proof}.}} Since $\iota_{g}:\mathrm{Mod}(\mathscr{O}_{X})\rightarrow
\mathrm{Mod}(\mathscr{A})$ is fully faithful,
\begin{itemize}
  \item [] $\mathrm{Hom}(\iota_{g}k(x)\otimes\Omega_{X}^{\mathscr{A}},\iota_{g}k(x))$
  \item [] $\simeq\mathrm{Hom}(\mathrm{gr}_{\hbar}(\iota_{g}k(x)\otimes\Omega_{X}^{\mathscr{A}}),k(x))$
  \item [] $\simeq\mathrm{Hom}(\mathrm{gr}_{\hbar}(\iota_{g}k(x))\otimes
\mathrm{gr}_{\hbar}\Omega_{X}^{\mathscr{A}},k(x))$
  \item [] $\simeq
\mathrm{Hom}((k(x)[1]\oplus{k}(x))\otimes\Omega_{X},k(x))$
  \item [] $\simeq
\mathrm{Hom}(k(x)[1]\oplus{k}(x),k(x))$
  \item [] $\simeq
\mathrm{Hom}(k(x)[1],k(x))\oplus\mathrm{Hom}(k(x),k(x))$
  \item [] $\simeq
\mathrm{Hom}(k(x),k(x))$.
\end{itemize}
Hence the identity of $k(x)$ defines a morphism
\begin{center}
$\varphi:\iota_{g}k(x)\otimes\Omega_{X}^{\mathscr{A}}\rightarrow\iota_{g}k(x).$
\end{center}
Since $\mathrm{gr}_{\hbar}(\iota_{g}k(x)\otimes\Omega_{X}^{\mathscr{A}})
\simeq\mathrm{gr}_{\hbar}(\iota_{g}k(x))$, $\mathrm{gr}_{\hbar}(\varphi)$
is an isomorphism which implies that $\varphi$ is an isomorphism
by Corollary 2.7. Hence we obtain
$\iota_{g}k(x)\simeq\iota_{g}k(x)\otimes\Omega_{X}^{\mathscr{A}}\simeq{S}_{X}
(\iota_{g}k(x))[-\mathrm{dim}X]$.$\hfill \square$\\

\noindent
{\bf{Proposition 5.8.}} \emph{Let $X$ and $Y$ be projective varieties. Assume
$F:\mathrm{D}^{\mathrm{b}}_{\mathrm{coh}}(\mathscr{A}_{X})
\stackrel{\sim}\longrightarrow\mathrm{D}^{\mathrm{b}}_{\mathrm{coh}}
(\mathscr{A}_{Y})$ is an equivalence, then $\mathrm{dim}X=\mathrm{dim}Y$.}\\

\noindent
{\emph{{Proof}.}}
First note that $F(\iota_{g}k(x))\simeq{F}(\iota_{g}k(x)\otimes\Omega_{X}^{\mathscr{A}})
=F(S_{X}(\iota_{g}k(x))[-\mathrm{dim}X])=S_{Y}F(\iota_{g}k(x))$
$[-\mathrm{dim}X]=F(\iota_{g}k(x))\otimes\Omega_{Y}^{\mathscr{A}}[\mathrm{dim}Y-\mathrm{dim}X].$

Taking $\mathrm{gr}_{\hbar}$ functor, we obtain $\mathrm{gr}_{\hbar}F(\iota_{g}k(x))
\simeq\mathrm{gr}_{\hbar}F(\iota_{g}k(x))\otimes\Omega_{Y}
[\mathrm{dim}Y-\mathrm{dim}X]$ which implies $\mathrm{dim}X=
\mathrm{dim}Y$ since $\mathrm{gr}_{\hbar}F(\iota_{g}k(x))\neq{0}$.$\hfill \square$\\

\noindent
{\bf{Definition 5.9.}} Let $\mathcal{D}$ be a $\mathbb{C}^{\hbar}$-linear
triangulated category with a Serre functor $S$. An object $P\in
\mathcal{D}$ is called point like of codimension $d$ if
\begin{itemize}
  \item [(i)] $S(P)\simeq{S}[d]$
  \item [(ii)] $\mathrm{Hom}(P,P[i])=0$ for $i<0$ and
  \item [(iii)] $k(P):=\mathrm{Hom}(P,P)\simeq\mathbb{C}$.
\end{itemize}
An object $P$ satisfying (iii) is called simple.\\

\noindent
{\bf{Example 5.10.}} For $x\in{X}$, $\iota_{g}k(x)$ is a point like of
codimension $\mathrm{dim}\mathrm{X}=n$ in $\mathrm{D}^{\mathrm{b}}_{\mathrm{coh}}
(\mathscr{A}_{X})$.
\begin{itemize}
  \item [(i)] $S(\iota_{g}k(x))\simeq\iota_{g}k(x)[n]$,
  \item [(ii)] $\mathrm{Hom}(\iota_{g}k(x),\iota_{g}k(x)[i])=0$ for $i<0$,
  \item [(iii)] $k(\iota_{g}k(x)):=\mathrm{Hom}(\iota_{g}k(x),\iota_{g}
k(x))\simeq\mathbb{C}$ by Propostion 5.6.
\end{itemize}

\noindent
{\bf{Lemma 5.11.}} \emph{Let $X$ be a topological space with sheaf of rings
$\mathscr{R}_{X}$. Suppose $\mathscr{M}\in\mathrm{D}^{\mathrm{b}}_{\mathrm{coh}}
(\mathscr{R}_{X})$ and $\mathrm{supp}(\mathscr{M})=Z_{1}\cup{Z}_{2}$ where
$Z_{1},Z_{2}\subset{X}$ are disjoint closed subsets. Then
$\mathscr{M}=\mathscr{M}_{1}\oplus\mathscr{M}_{2}$ with
$\mathrm{supp}\mathscr{M}_{j}\subset{Z}_{j}$ for $j=1,2$. In particular,
the same result holds when replacing the sheaf of rings $\mathscr{R}_{X}$
by the DQ-algebroid stack $\mathscr{A}_{X}$.}\\

\noindent
{\emph{{Proof}.}} We have the following distinguished triangle
\begin{center}
$R\Gamma_{Z_{1}\cap{Z}_{2}}(\mathscr{M})\longrightarrow
{R}\Gamma_{Z_{1}}(\mathscr{M})\oplus{R}\Gamma_{Z_{2}}(\mathscr{M})\longrightarrow
{R}\Gamma_{Z_{1}\cup{Z}_{2}}(\mathscr{M})\stackrel{+1}\longrightarrow$.
\end{center}
Since $R\Gamma_{Z_{1}\cap{Z}_{2}}(\mathscr{M})=0$ and
$R\Gamma_{Z_{1}\cup{Z}_{2}}(\mathscr{M})=\mathscr{M}$, we get
\begin{center}
$\mathscr{M}\simeq{R}\Gamma_{Z_{1}}(\mathscr{M})\oplus{R}\Gamma_{Z_{2}}
(\mathscr{M})$.
\end{center}
It remains to prove that $R\Gamma_{Z_{1}}(\mathscr{M})$ and
$R\Gamma_{Z_{2}}(\mathscr{M})$ are coherent. We have the following
two exact sequences
\begin{center}
$0\longrightarrow{H}^{i}_{Z_{1}}(\mathscr{M})\longrightarrow{H}^{i}
(\mathscr{M})\longrightarrow{H}^{i}_{Z_{2}}(\mathscr{M})\longrightarrow{0}$
\end{center}
 \begin{center}
 $0\longrightarrow{H}^{i}_{Z_{2}}(\mathscr{M})\longrightarrow{H}^{i}
 (\mathscr{M})\longrightarrow{H}^{i}_{Z_{1}}(\mathscr{M})\longrightarrow{0}$.
 \end{center}
Since $H^{i}(\mathscr{M})$ is coherent, we obtain $H^{i}_{Z_{1}}(\mathscr{M})$
and $H^{i}_{Z_{2}}(\mathscr{M})$ are coherent.  $\hfill \square$\\

\noindent
{\bf{Lemma 5.12.}} \emph{Suppose $\mathscr{M}$ is a simple object in
$\mathrm{D}^{\mathrm{b}}_{\mathrm{coh}}(\mathscr{A}_{X})$ with
$\hbar{H}^{i}\mathscr{M}=0$ and zero dimensional support.
If $\mathrm{Hom}_{\mathscr{A}_{X}}(\mathscr{M},\mathscr{M}[i])=0$
for $i<0$, then $\mathscr{M}\simeq\iota_{g}k(x)[m]$ for some
closed point $x\in{X}$ and some integer $m$.}\\

\noindent
{\emph{{Proof}.}} Let us show that $\mathscr{M}$ is concentrated in exactly
one closed point. Otherwise $\mathscr{M}$ could be written as direct sum
$\mathscr{M}\simeq\mathscr{M}_{1}\oplus\mathscr{M}_{2}$ with
$\mathscr{M}_{j}\simeq{0},j=1,2$ by Lemma 5.11. However, a direct sum of
two non-trivial complexes is never simple, Indeed, the projection to one of
the two summands is not invertible.

Thus, we may assume that the support of all cohomology sheaves
$\mathscr{H}^{i}$ of $\mathscr{M}$ consists of the same closed
point $x\in{X}$. Set
\begin{center}
$m_{0}:=\mathrm{max}\{i|\mathscr{H}^{i}\neq{0}\}$ and $m_{1}:=
\mathrm{min}\{i|\mathscr{H}^{i}\neq{0}\}$.
\end{center}
Since both modules $\mathscr{H}^{m_{0}}=\iota_{g}\mathscr{F}_{0}$ and
$\mathscr{H}^{m_{1}}=\iota_{g}\mathscr{F}_{1}$ are concentrated in
$x\in{X}$ and there exists non-trivial homomorphism $\mathscr{F}_{0}
\rightarrow\mathscr{F}_{1}$,  there exists non-trivial
homomoprhism $\mathscr{H}^{m_{0}}\longrightarrow
\mathscr{H}^{m_{1}}$ because $\iota_{g}$ is fully faithful. Now the
composition
\begin{center}
$\mathscr{M}[m_{0}]\longrightarrow\mathscr{H}^{m_{0}}\longrightarrow
\mathscr{H}^{m_{1}}\longrightarrow\mathscr{M}[m_{1}]$
\end{center}
is nontrivial. This shows $m_{0}=m_{1}$ by Definition 5.9 (ii).
Hence $\mathscr{M}=\mathscr{N}[m]$ with $\mathscr{N}\simeq\iota_{g}\mathscr{F}$
for some coherent $\mathscr{O}_{X}$-module $\mathscr{F}$ supported in
$x$ and $m=m_{0}=m_{1}$. Since
$\mathrm{Hom}_{\mathscr{A}_{X}}(\iota_{g}\mathscr{F},\iota_{g}
\mathscr{F})\simeq\mathrm{Hom}_{\mathscr{O}_{X}}
(\mathrm{gr}_{\hbar}\iota_{g}\mathscr{F},\mathscr{F})\simeq
\mathrm{Hom}_{\mathscr{O}_{X}}(\mathscr{F},\mathscr{F})\simeq\mathbb{C}$,
$\mathscr{F}\simeq{k}(x)$ by Proposition 4.4.$\hfill \square$\\

\noindent
{\bf{Proposition 5.13.}} \emph{Let $X$ be a projective variety. Suppose that
$\Omega_{X}$ or $\Omega_{X}^{*}$ is ample. Then the point like objects
in $\mathrm{D}^{\mathrm{b}}_{\mathrm{coh}}(\mathscr{A}_{X})$ are the objects
which are isomorphic to $\iota_{g}k(x)[m]$, where $x\in{X}$ is a closed
point and $m\in\mathbb{Z}$.}\\

\noindent
{\emph{{Proof}.}} Let $\mathscr{P}$ be a point like object in
$\mathrm{D}^{\mathrm{b}}_{\mathrm{coh}}(\mathscr{A}_{X})$. Then
$\hbar:\mathscr{P}\longrightarrow\mathscr{P}$ is zero.
Indeed, if $\hbar:\mathscr{P}\longrightarrow\mathscr{P}$ is not zero
then it is an isomorphism since $\mathrm{Hom}(\mathscr{P},\mathscr{P})$
is a field. Then from the distinguished triangle
\begin{center}
 $\mathscr{P}\stackrel{\hbar}\longrightarrow\mathscr{P}
\longrightarrow\mathrm{gr}_{\hbar}\mathscr{P}\stackrel{+1}\longrightarrow$
\end{center}
 we get $\mathrm{gr}_{\hbar}\mathscr{P}\simeq0$ which implies that
$\mathscr{P}\simeq{0}$ a contradiction. Hence we get
$\hbar{H}^{j}\mathscr{P}=0$ for each $j$ and $H^{j}(\mathscr{P})\simeq
\iota_{g}\mathscr{F}$ for some coherent $\mathscr{O}_{X}$-module $\mathscr{F}$.
Since $\mathscr{P}\otimes\Omega_{X}^{\mathscr{A}}\simeq\mathscr{P},
\mathrm{gr}_{\hbar}\mathscr{P}\otimes\Omega_{X}\simeq
\mathrm{gr}_{\hbar}\mathscr{P}$. Hence $\mathrm{dim}\mathrm{Supp}\mathscr{P}
=0$ since $\mathrm{supp}\mathscr{P}=\mathrm{supp}\mathrm{gr}_{\hbar}\mathscr{P}$.
Hence by Lemma 5.12, we obtain $\mathscr{P}\simeq\iota_{g}k(x)[m]$.   $\hfill \square$\\

\noindent
{\bf{Lemma 5.14.}} \emph{Suppose that $F:\mathrm{D}^{\mathrm{b}}_{\mathrm{coh}}
(\mathscr{A}_{X})\longrightarrow\mathrm{D}^{\mathrm{b}}_{\mathrm{coh}}
(\mathscr{A}_{Y})$ is an exact equivalence with $\Omega_{X}$ an
(anti)-ample line bundle. Then there exists a bijection between
point like objects on $X$ and point like objects on $Y$ and invertible
modules on $X$ and invertible modules on $Y$ up to a shift functor.} \\

\noindent
{\emph{{Proof}.}} First, it is easy to check that there is a bijection
between point like objects and point like objects.

Next we shall prove that point like objects in $\mathrm{D}^{\mathrm{b}}
_{\mathrm{coh}}(\mathscr{A}_{Y})$ are of the form $\iota_{g}k(y)[m]$
for some $y\in{Y}$ and $m\in\mathbb{Z}$.

Since $F$ is an equivalence, for any closed point $y\in{Y}$ there exists
a closed point $x_{y}\in{X}$ and an integer $m_{y}$ such that
$\iota_{g}k(y)\simeq{F}(\iota_{g}k(x_{y})[m_{y}])$.

Suppose there exists a point like object $P\in\mathrm{D}^{\mathrm{b}}
_{\mathrm{coh}}(\mathscr{A}_{Y})$ which is not of the form $\iota_{g}
k(y)[m]$ and we denote by $x_{P}\in{X}$ the closed point such that
$F(\iota_{g}k(x_{P})[m_{P}])\simeq{P}$ for a certain
$m_{P}\in\mathbb{Z}$.
Note that $x_{P}\neq{x}_{y}$ for all $y\in{Y}$. Hence we have for all
$y\in{Y}$ and all $m\in\mathbb{Z}$
\begin{itemize}
  \item [] $\mathrm{Hom}(P,\iota_{g}k(y)[m])$
  \item [] $\simeq\mathrm{Hom}(F(\iota_{g}k(x_{P}))[m_{P}],
F(\iota_{g}k(x_{y}))[m_{y}+m])$
  \item [] $\simeq\mathrm{Hom}(\iota_{g}k(x_{P}),\iota_{g}k(x_{y})
[m_{y}+m-m_{P}])$
  \item [] $\simeq{0}$.
\end{itemize}
Since the objects $\iota_{g}k(y)[m]$ form a spanning class in
$\mathrm{D}^{\mathrm{b}}_{\mathrm{coh}}(\mathscr{A}_{Y})$, this
implies $P\simeq{0}$ which is absurd. Hence, point like objects
in $\mathrm{D}^{\mathrm{b}}_{\mathrm{coh}}(\mathscr{A}_{Y})$
are exactly the objects of the form $\iota_{g}k(y)[m]$.

Now we prove that if $\mathscr{P}$ is an invertible $\mathscr{A}_{X}$-module
then $F(\mathscr{P})$ is an invertible $\mathscr{A}_{Y}$-module up
to a shift functor.

Since
\begin{itemize}
  \item [] $\mathrm{Hom}_{\mathscr{O}_{X}}(\mathrm{gr}_{\hbar}\mathscr{P},
  k(x)[j])$
  \item [] $\simeq\mathrm{Hom}_{\mathscr{A}_{X}}
  (\mathscr{P},\iota_{g}k(x)[j])$
  \item [] $\simeq\mathrm{Hom}_{\mathscr{A}_{Y}}(F(\mathscr{P}),F(\iota_{g}k(x))[j])$
  \item [] $\simeq\mathrm{Hom}_{\mathscr{A}_{Y}}(F(\mathscr{P}),\iota_{g}k(y)[j+m])$
  \item [] $\simeq\mathrm{Hom}_{\mathscr{O}_{Y}}(\mathrm{gr}_{\hbar}F(\mathscr{P}),k(y)[j+m])$
  \item [] $\simeq\mathrm{Hom}_{\mathscr{O}_{Y}}(H^{-m}(\mathrm{gr}_{\hbar}
  (F(\mathscr{P}))),k(y)[j])$,
\end{itemize}
this implies that $\mathrm{gr}_{\hbar}F(\mathscr{P})\simeq{L}[m]$ for
some line bundle $L$ by Corollary 2.19, Proposition 4.5, and Proposition
4.7. Hence $F(\mathscr{P})[-m]$ is invertible, we get the result.      $\hfill \square$\\

\noindent
{\bf{Proposition 5.15.}} \emph{Let $X$ and $Y$ be two projective varieties. Assume
$F:\mathrm{D}^{\mathrm{b}}_{\mathrm{coh}}(\mathscr{O}_{X}^{\hbar})\longrightarrow
\mathrm{D}^{\mathrm{b}}_{\mathrm{coh}}(\mathscr{O}_{Y}^{\hbar})$ is fully
faithful, then it induces a fully faithful functor $F':=
\mathrm{gr}_{\hbar}\circ{F}\circ(\cdot)^{\mathrm{R}\hbar}:
\mathrm{D}^{\mathrm{b}}_{\mathrm{coh}}(\mathscr{O}_{X})\longrightarrow
\mathrm{D}^{\mathrm{b}}_{\mathrm{coh}}(\mathscr{O}_{Y})$.}\\

\noindent
{\emph{{Proof}.}}  We have
\begin{center}
    \begin{itemize}
      \item [] $\mathrm{RHom}_{\mathscr{O}_{Y}}
(\mathrm{gr}_{\hbar}F(\mathscr{M})^{\mathrm{R}\hbar},\mathrm{gr}_{\hbar}F(\mathscr{N})^{\mathrm{R}\hbar})$
      \item []  $\simeq\mathbb{C}\stackrel{\mathrm{L}}\otimes_{\mathbb{C}^{\hbar}}
\mathrm{RHom}_{\mathscr{A}_{Y}}(F(\mathscr{M})^{\mathrm{R}\hbar},F(\mathscr{N})^{\mathrm{R}\hbar})$
      \item []   $\simeq\mathbb{C}\stackrel{\mathrm{L}}\otimes_{\mathbb{C}^{\hbar}}
\mathrm{RHom}_{\mathscr{A}_{X}}((\mathscr{M})^{\mathrm{R}\hbar},(\mathscr{N})^{\mathrm{R}\hbar})$
      \item []  $\simeq\mathrm{RHom}_{\mathscr{O}_{X}}(\mathrm{gr}_{\hbar}(\mathscr{M})^{\mathrm{R}\hbar},
\mathrm{gr}_{\hbar}(\mathscr{N})^{\mathrm{R}\hbar})$
      \item []   $\simeq\mathrm{RHom}_{\mathscr{O}_{X}}(\mathscr{M},\mathscr{N})$.$\hfill \square$\\
    \end{itemize}
\end{center}

Now we are ready to prove Theorem 1.2.\\

 \noindent
 {\bf{Theorem 5.16.}}\emph{Let $X$ and $Y$ be two
 smooth projective varieties and assume that the (anti)-canonical
 bundle of $X$ is ample. If there exists an exact equivalence
 $\mathrm{D}^{\mathrm{b}}_{\mathrm{coh}}(\mathscr{O}^{\hbar}_{Y})
 \stackrel{\sim}\longrightarrow
 \mathrm{D}^{\mathrm{b}}_{\mathrm{coh}}(\mathscr{O}_{X}^{\hbar})$, then
 $X\simeq{Y}$. In particular, the (anti)-canonical bundle of
 $Y$ is also ample}. \\

\noindent
{\emph{{Proof}.}} We first assume that $\Omega_{Y}$ is ample. By Corollary 2.18
and Lemma 5.14, we may assume that $F(\mathscr{O}_{X}^{\hbar})\simeq
\mathscr{O}_{Y}^{\hbar}$. Since $F$ is an equivalence, $\mathrm{dim}X=
\mathrm{dim}Y=n$ by Proposition 5.8, this proves
\begin{center}
$F((\Omega_{X}^{\hbar})^{\otimes{k}})\simeq
{F}(S_{X}^{k}(\mathscr{O}_{X}^{\hbar}))[-kn]
\simeq{S}_{Y}^{k}(F(\mathscr{O}_{X}^{\hbar}))[-kn]\simeq
{S}_{Y}^{k}(\mathscr{O}_{Y}^{\hbar})[-kn]\simeq
(\Omega_{Y}^{\hbar})^{\otimes{k}}.$
\end{center}
Using the fact that $F'$ is fully faithful by Theorem 5.15, we have
\begin{center}
$H^{0}(X,\Omega_{X}^{\otimes{k}})=
\mathrm{Hom}(\mathscr{O}_{X},\Omega_{X}^{\otimes{k}})
\simeq\mathrm{Hom}(F'(\mathscr{O}_{X}),F'(\Omega_{X}^{\otimes{k}}))\simeq
\mathrm{Hom}_{\mathscr{O}_{Y}}(\mathscr{O}_{Y},\Omega_{Y}^{\otimes{k}})
\simeq{H}^{0}(Y,\Omega_{Y}^{\otimes{k}})$.
\end{center}

The product in $\bigoplus{H}^{0}(X,\Omega_{X}^{\otimes{k}})$ can be
expressed in terms of compositions:
namely, for $s_{i}\in{H}^{0}(X,\Omega_{X}^{\otimes{k}_{i}})\simeq
\mathrm{Hom}(\mathscr{O}_{X},\Omega_{X}^{\otimes{k}_{i}})$ one has
\begin{center}
$s_{1}\cdot{s}_{2}=S_{X}^{k_{1}}(s_{2})[-k_{1}n]\circ{s}_{1}$
\end{center}
and similarly for sections on $Y$.

Hence, the induced bijection
\begin{center}
$\bigoplus{H}^{0}(X,\Omega_{X}^{\otimes{k}})\simeq\bigoplus{H}^{0}(Y,\Omega_{Y}^{\otimes{k}})$
\end{center}
is a ring homomorphism. Since $Y$ is also ample, then it shows
\begin{center}
$X\simeq\mathrm{Proj}(\bigoplus{H}^{0}(X,\Omega_{X}^{\otimes{k}}))\simeq
\mathrm{Proj}(\bigoplus{H}^{0}(Y,\Omega_{Y}^{\otimes{k}}))\simeq{Y}$.
\end{center}
Thus, under the assumption that $\Omega_{Y}$ (or $\Omega_{Y}^{*}$) is ample,
we have proved the assertion.

Finally, we shall remove the assumption of the ampleness of $\Omega_{Y}$ by proving
that it is really an ample line bundle.

We continue to use that for any $k(y)$ with $y\in{Y}$ a closed point, there
exists a closed point $x_{y}\in{X}$ such that $F(\iota_{g}k(x_{y}))=
\iota_{g}k(y)$ (indeed, using $\mathrm{Hom}(\mathscr{O}_{X},k(x_{y}))\simeq
\mathrm{Hom}(\mathscr{O}_{X}^{\hbar},\iota_{g}k(x_{y}))\simeq
\mathrm{Hom}(\mathscr{O}_{Y}^{\hbar},\iota_{g}k(y)[m])\simeq
\mathrm{Hom}(\mathscr{O}_{Y},k(y)[m])$ would implies that $m=0$) and
that $F((\Omega_{X}^{\hbar})^{\otimes{k}})\simeq(\Omega_{Y}^{\hbar})^{\otimes{k}}$
for all $k\in\mathbb{Z}$.

Now we show that $F'(k(x_{y}))\simeq{k}(y)$ and that
$F'(\Omega_{X}^{\otimes{k}})\simeq\Omega_{Y}^{\otimes{k}}$ for all $k\in\mathbb{Z}$.
Indeed, for any $k(y)$ with $y\in{Y}$ a closed point, we have
\begin{itemize}
  \item [] $\mathrm{Hom}(\mathrm{gr}_{\hbar}Fk(x_{y})[[\hbar]],k(y)[i])$
  \item [] $\simeq\mathrm{Hom}(Fk(x_{y})[[\hbar]],\iota_{g}k(y)[i])$
  \item [] $\simeq\mathrm{Hom}(k(x_{y})[[\hbar]],\iota_{g}k(x_{y})[i])$
  \item [] $\simeq\mathrm{Hom}(k(x_{y}),k(x_{y})[i])$.
\end{itemize}
Hence by Lemma 4.8, $\mathscr{F}:=\mathrm{gr}_{\hbar}Fk(x_{y})[[\hbar]]\simeq
{F}'k(x_{y})$ is concentrated in $y\in{Y}$. Hence
$\mathrm{gr}_{\hbar}Fk(x_{y})[[\hbar]]$ is a simple object in $\mathrm{D}^{\mathrm{b}}
_{\mathrm{coh}}(\mathscr{O}_{X})$ with zero dimension support and
$\mathrm{Hom}(\mathscr{F},\mathscr{F}[i])\simeq{0}$ for $i<0$.
Then  $\mathscr{F}\simeq{k}(y)$
by Proposition 4.4.

Also, we get
\begin{center}
            $\mathrm{Hom}(k(x_{y}),k(x_{y})[1])\simeq
\mathrm{Hom}(F'k(x_{y}),F'k(x_{y})[1])\simeq
\mathrm{Hom}(k(y),k(y)[1])$.
\end{center}
Now the proof for ampleness of $\Omega_{Y}$ is the same as [Huy06] P.97 and
P.98.    $\hfill \square$\\

\noindent
{\bf{Corollary 5.17.}} \emph{Let $C$ be a curve of genus $g(C)\neq{1}$ and let
$Y$ be a projective variety. If there exists an exact equivalence
$\mathrm{D}^{\mathrm{b}}_{\mathrm{coh}}(\mathscr{O}_{C}^{\hbar})
\stackrel{\sim}\longrightarrow\mathrm{D}^{\mathrm{b}}_{\mathrm{coh}}(\mathscr{O}_{Y}^{\hbar})$,
then $Y$ is a curve isomorphic to $C$.}

 \section{Proof of Theorem 1.4}
 \label{sec.recall}
We first review the notions of generators and compactness in triangulated
categories.\\

\noindent {\bf{Definition 6.1.}} Let $\mathcal{T}$ be a
triangulated category. Let $\mathcal{I}=(C_{i})_{i\in{I}}$ be a set
of objects of $\mathcal{T}$. One says that $\mathcal{I}$ generates
$\mathcal{T}$ if for every $D\in\mathcal{T}$ such that
$\mathrm{Hom}_{\mathcal{T}}(C_{i}[n],D)=0$ for every
$C_{i}\in\mathcal{I}$ and $n\in\mathbb{Z}$ then $D\simeq{0}$. \\

\noindent {\bf{Definition 6.2.}} Let $\mathcal{T}$ be a
triangulated category. Assume that $\mathcal{T}$ admits arbitrary
small direct sums. An object $T$ in $\mathcal{T}$ is compact if
$\mathrm{Hom}_{\mathcal{T}}(T,\cdot)$ commutes with small direct
sums. We denote by $\mathcal{T}^{\mathrm{c}}$ the full subcategory
of $\mathcal{T}$ consisting of compact objects.\\

Let $X$ be a variety, we denote by $\mathrm{D}_{\mathrm{qcoh}}
(\mathscr{O}_{X})$ the derived category of $\mathscr{O}_{X}$-modules
with quasi-coherent cohomology sheaves and $\mathrm{Perf}(X)$
the category of perfect complex on $X$. Since $X$ is smooth,
we have $\mathrm{Perf}(X)=\mathrm{D}^{\mathrm{b}}_{\mathrm{coh}}
(\mathscr{O}_{X})$.\\

\noindent {\bf{Example 6.3.}} ([BvdB03, Theorem 3.3.1]).
For a variety $X$, we have $\mathrm{D}_{\mathrm{qcoh}}(\mathscr{O}_{X})^{\mathrm{c}}
=\mathrm{Perf}(X)=\mathrm{D}^{\mathrm{b}}_{\mathrm{coh}}
(\mathscr{O}_{X})$.\\

\noindent {\bf{Definition 6.4.}} A triangulated category $\mathcal{T}$
is called compactly generated if it admits arbitrary small direct sums
and is generated by a set $\mathcal{I}$ of compact objects of $\mathcal{T}$.\\

\noindent {\bf{Theorem 6.5.}} (The Brown Representability Theorem
[Nee96, Theorem 4.1]). \emph{Let $\mathcal{S}$ be a compactly generated
triangulated category, $\mathcal{T}$ any  triangulated category. Let
$F:\mathcal{S}\rightarrow\mathcal{T}$ be a triangulated functor.
Suppose $F$ commutes with small direct sums; that is,
the natural map
\begin{center}
$\bigoplus\limits_{\lambda\in\Lambda}
{F}(s_{\lambda})\rightarrow{F}(\bigoplus\limits_{\lambda\in\Lambda}
s_{\lambda})$
 \end{center}
is an isomorphism. Note that although we are not assuming that
$\mathcal{T}$ has direct sums, we are assuming that the object on
the right is a direct sum of the objects on the left. Then $F$ has a
right adjoint $G:\mathcal{T}\rightarrow\mathcal{S}$.}\\

\noindent {\bf{Example 6.6.}} ([BvdB03, Theorem 3.3.1]).
$\mathrm{D}_{\mathrm{qcoh}}(\mathscr{O}_{X})$ is compactly generated
by an object $\mathscr{M}\in\mathrm{D}^{\mathrm{b}}_{\mathrm{coh}}
(\mathscr{O}_{X})$.\\

Now we begin to prove Theorem 1.4 assuming the following equivalence
\begin{center}
(6.1)\hfill$F:\mathrm{D}^{\mathrm{b}}_{\mathrm{coh}}
(\mathscr{O}_{X})\stackrel{\sim}\longrightarrow
\mathrm{D}^{\mathrm{b}}_{\mathrm{coh}}(\mathscr{O}_{Y}).\hfill$
 \end{center}
For a scheme $X$, we denote by $\mathrm{D}(Qcoh(X))$ the derived category of
quasi-coherent sheaves on $X$, $\mathrm{D}^{\mathrm{b}}
(coh(X))$ the bounded derived category of coherent sheaves on $X$ and
$\mathrm{Perf}(X)$ the perfect complex on $X$.\\

Recall that for schemes $X,Y$ and $\mathscr{E}\in
\mathrm{D}(Qcoh(X\times{Y}))$, we have the following
Fourier-Mukai transform:
\begin{center}
$\Phi_{\mathscr{E}}(-):={\mathrm{R}p_{2}}_{*}
(\mathscr{E}\stackrel{\mathrm{L}}\otimes{p}_{1}^{*}(-)):
\mathrm{D}(Qcoh(X))\rightarrow\mathrm{D}(Qcoh(Y))$
 \end{center}
where $p_{1}:X\times{Y}\rightarrow{X}$
and $p_{2}:X\times{Y}\rightarrow{Y}$  are projections.\\

We have the following results.\\

\noindent {\bf{Theorem 6.7.}} ([LO10, Corollary 9.13]). \emph{Let $X$ be a
quasi-projective scheme and $Y$ be quasi-compact and separated
scheme. Let $F:\mathrm{Perf}(X)\rightarrow\mathrm{D}(Qcoh(Y))$
(we keep the same notation as (6.1)) be a
fully faithful functor. Then there is an object
$\mathscr{E}\in\mathrm{D} (Qcoh(X\times{Y}))$ such that
\begin{itemize}
  \item [(1)] the functor $\Phi_{\mathscr{E}}|_{\mathrm{Perf}X}:
\mathrm{Perf}(X)\rightarrow\mathrm{D}(Qcoh{Y})$ is fully faithful
and $\Phi_{\mathscr{E}}(P)\simeq{F}(P)$ for any $P\in\mathrm{Perf}(X)$,
  \item [(2)] if $F$ sends $\mathrm{Perf}(X)$ to $\mathrm{Perf}(Y)$ then
the functor $\Phi_{\mathscr{E}}:\mathrm{D}(Qcoh{X})\rightarrow
\mathrm{D}(Qcoh{Y})$ is fully faithful and also sends
$\mathrm{Perf}(X)$ to $\mathrm{Perf}(Y)$,
  \item [(3)] if $Y$ is noetherian and $F$ sends $\mathrm{Perf}(X)$
to $\mathrm{D}^{\mathrm{b}}(Qcoh{Y})_{\mathrm{coh}}$, then the
object $\mathscr{E}$ is isomorphic to an object of
$\mathrm{D}^{\mathrm{b}} (coh(X\times{Y}))$.
 \end{itemize}}

Now we assume that $X,Y$ are smooth varieties. Recall that we have
the following equivalences: \begin{center}
$D(Qcoh{X})\stackrel{\sim}\longrightarrow\mathrm{D}_{\mathrm{qcoh}}
(\mathscr{O}_{X})$ and $\mathrm{D}^{\mathrm{b}}(coh{X})
\stackrel{\sim}\longrightarrow\mathrm{D}^{\mathrm{b}}
_{\mathrm{coh}}(\mathscr{O}_{X})$.
 \end{center}

\noindent {\bf{Corollary 6.8.}} \emph{By Theorem 6.7, there exists
$\mathscr{E}\in\mathrm{D}^{\mathrm{b}}_{\mathrm{coh}}
(\mathscr{O}_{X\times{Y}})$ such that
\begin{center}
$\Phi_{\mathscr{E}}:\mathrm{D}_{\mathrm{qcoh}}(\mathscr{O}_{X})
\longrightarrow{D}_{\mathrm{qcoh}}(\mathscr{O}_{Y})$
\end{center}
is fully faithful and $\Phi_{\mathscr{E}}(P)\simeq{F}(P)$ for
$P\in\mathrm{D}^{\mathrm{b}}_{\mathrm{coh}}(\mathscr{O}_{X})$.}\\

\noindent {\bf{Definition 6.9. }} We say $\mathscr{M}\in
\mathrm{D}(\mathscr{O}_{X}^{\hbar})$ is quasi-coherent if
$\mathrm{gr}_{\hbar}\mathscr{M}\in\mathrm{D}_{\mathrm{qcoh}}
(\mathscr{O}_{X})$.\\

We denote by $\mathrm{D}_{\mathrm{qcc}}(\mathscr{O}_{X}^{\hbar})$
the full subcategory of $\mathrm{D}(\mathscr{O}_{X}^{\hbar})$
consisting of cohomologically complete and quasi-coherent $\mathscr{O}_{X}^{\hbar}$-modules.\\

\noindent {\bf{Proposition 6.10.}} The category $\mathrm{D}_{\mathrm{qcc}}
(\mathscr{O}_{X}^{\hbar})$ is a $\mathbb{C}^{\hbar}$-linear full triangulated
subcategory of $\mathrm{D}(\mathscr{O}_{X}^{\hbar})$.\\

\noindent {\bf{Corollary 6.11.}} \emph{$\mathrm{D}_{\mathrm{qcc}}(\mathscr{O}_{X}^{\hbar})$
is compactly generated by $\mathscr{M}^{\mathrm{R}\hbar}$ for some
$\mathscr{M}\in\mathrm{D}^{\mathrm{b}}_{\mathrm{coh}}(\mathscr{O}_{X})$.}\\

\noindent
{\emph{{Proof}}.} Use Proposition 2.10 and Example 6.6.\\

\noindent {\bf{Proposition 6.12.}} ([KS12, Proposition 1.5.8]). \emph{Let $\mathscr{M}\in\mathrm{D}(\mathscr{O}_{X}^{\hbar})$
be a cohomologically complete object and $a\in\mathbb{Z}$. If
$H^{i}(\mathrm{gr}_{\hbar}(\mathscr{M}))=0$ for any $i<a$, then
$H^{i}(\mathscr{M})=0$ for any $i<a$.}\\

\noindent {\bf{Proposition 6.13.}} ([KS12, Theorem 1.6.4]). \emph{If $\mathscr{M}\in\mathrm{D}^{+}_{\mathrm{qcc}}
(\mathscr{O}_{X}^{\hbar})$ is such that $\mathrm{gr}_{\hbar}\mathscr{M}\in
\mathrm{D}^{+}_{\mathrm{coh}}(\mathscr{O}_{X})$, then
$\mathscr{M}\in\mathrm{D}^{+}_{\mathrm{coh}}(\mathscr{O}_{X}^{\hbar})$.}\\

\noindent {\bf{Definition 6.14. }} We denote by $(\cdot)^{cc}$ the functor
\begin{center}
$\mathrm{R}\mathscr{H}om_{\mathscr{O}_{X}^{\hbar}}
(\mathscr{O}_{X}^{\mathrm{loc}}/\mathscr{O}_{X}^{\hbar}[-1],\cdot):
\mathrm{D}(\mathscr{O}_{X}^{\hbar})\rightarrow\mathrm{D}(\mathscr{O}_{X}^{\hbar})$.
\end{center}
\noindent {\bf{Corollary 6.15.}} ([Pet11, Corollary 3.2.7]).
For every $\mathscr{M}\in\mathrm{D}(\mathscr{O}_{X}^{\hbar})$, we have
\begin{center}
$\mathrm{gr}_{\hbar}(\mathscr{M}^{cc})\simeq\mathrm{gr}_{\hbar}
\mathscr{M}$.
\end{center}
\noindent {\bf{Definition 6.16.}} Let $(\mathscr{M}_{i})_{i\in{I}}$
be a family of objects of $\mathrm{D}_{\mathrm{cc}}(\mathscr{O}_{X}^{\hbar})$. We set
\begin{center}
$\overline{\bigoplus\limits_{i\in{I}}}\mathscr{M}_{i}=
(\bigoplus\limits_{i\in{I}}\mathscr{M}_{i})^{cc}$
\end{center}
where $\bigoplus$ denotes the direct sum in the category $\mathrm{D}(\mathscr{O}_{X}^{\hbar})$.\\

\noindent {\bf{Proposition 6.17.}} ([Pet11, Proposition 3.2.9]). \emph{The category
$\mathrm{D}_{\mathrm{cc}}(\mathscr{O}_{X}^{\hbar})$ admits direct sums.
The direct sum of the family $(\mathscr{M}_{i})_{i\in{I}}$ is given
by $\overline{\bigoplus\limits_{i\in{I}}}\mathscr{M}_{i}$.}\\

\noindent {\bf{Proposition 6.18.}} ([Pet11, Proposition 3.2.10]).
\emph{The category $\mathrm{D}_{\mathrm{qcc}}(\mathscr{O}_{X}^{\hbar})$
admits direct sums. The direct sum of the family $(\mathscr{M}_{i})_{i\in{I}}$
is given by $\overline{\bigoplus\limits_{i\in{I}}}\mathscr{M}_{i}$.}\\

\noindent {\bf{Proposition 6.19.}} ([BvdB03, Theorem 3.3.3]). \emph{Let
$X$ be a quasi-compact, separated scheme and let $Y$ be a scheme.
Let $f:X\rightarrow{Y}$ be a separated morphism.
Let $Rf_{*}:\mathrm{D}_{\mathrm{qcoh}}(\mathscr{O}_{X})
\rightarrow\mathrm{D}_{\mathrm{qcoh}}(\mathscr{O}_{Y})$ be the direct
image functor. Then the natural map
\begin{center}
$\bigoplus\limits_{i\in{I}}Rf_{*}(\mathscr{M}_{i})\rightarrow
{R}f_{*}(\bigoplus\limits_{i\in{I}}\mathscr{M}_{i})$
\end{center}
is an isomorphism; that is, $Rf_{*}$ commutes with direct sums.}\\

\noindent {\bf{Corollary 6.20.}} ([KS12, Proposition 1.5.12]). \emph{Let $f:X\rightarrow{Y}$
be a morphism of varieties and let $\mathscr{M}\in
\mathrm{D}(f^{-1}\mathscr{O}_{Y}^{\hbar})$, then}
\begin{center}
$Rf_{*}\mathscr{M}^{cc}\simeq(Rf_{*}\mathscr{M})^{cc}$ in $\mathrm{D}
(\mathscr{O}_{Y}^{\hbar})$.
\end{center}
Define a functor
\begin{center}
$\Phi_{\mathscr{E}^{\mathrm{R}\hbar}}:
\mathrm{D}_{\mathrm{qcc}}(\mathscr{O}_{X}^{\hbar})\rightarrow
\mathrm{D}_{\mathrm{qcc}}(\mathscr{O}_{Y}^{\hbar})$
\end{center}
by $\Phi_{\mathscr{E}^{\mathrm{R}\hbar}}
(\mathscr{M}):={\mathrm{R}p_{2}}_{*}(\mathscr{E}^{\mathrm{R}\hbar}
\stackrel{\mathrm{L}}\otimes_{p_{1}^{-1}\mathscr{O}_{X}^{\hbar}}p_{1}^{-1}
\mathscr{M})^{cc}$ where $p_{1}:X\times{Y}\rightarrow{X}$
and $p_{2}:X\times{Y}\rightarrow{Y}$  are projections.\\

\noindent {\bf{Proposition 6.21.}} \emph{The functor $\Phi_{\mathscr{E}^{\mathrm{R}\hbar}}$
commutes with direct sums, i.e., the following natural morphism is an isomorphism:}
\begin{itemize}
  \item [] $(\bigoplus{\mathrm{R}p_{2}}_{*}(\mathscr{E}^{\mathrm{R}\hbar}
  \stackrel{\mathrm{L}}\otimes_{p_{1}^{-1}\mathscr{O}_{X}^{\hbar}}p_{1}^{-1}\mathscr{M}_{i}))^{cc}$
  \item [] $\longrightarrow({\mathrm{R}p_{2}}_{*}(\bigoplus(\mathscr{E}^{\mathrm{R}\hbar}
  \stackrel{\mathrm{L}}\otimes_{p_{1}^{-1}\mathscr{O}_{X}^{\hbar}}p_{1}^{-1}\mathscr{M}_{i})))^{cc}$
  \item [] $\stackrel{\sim}\longrightarrow{\mathrm{R}p_{2}}_{*}(\bigoplus(\mathscr{E}^{\mathrm{R}\hbar}
  \stackrel{\mathrm{L}}\otimes_{p_{1}^{-1}\mathscr{O}_{X}^{\hbar}}p_{1}^{-1}\mathscr{M}_{i}))^{cc}$
  \item [] $\longrightarrow{\mathrm{R}p_{2}}_{*}(\mathscr{E}^{\mathrm{R}\hbar}
  \stackrel{\mathrm{L}}\otimes_{p_{1}^{-1}\mathscr{O}_{X}^{\hbar}}p_{1}^{-1}
  (\bigoplus\mathscr{M}_{i})^{cc})^{cc}.$
\end{itemize}
\emph{{Proof}}. The first arrow is the natural morphism. The second isomorphism
follows from Corollary 6.20. Finally we have the following morphisms
\begin{itemize}
  \item [] ${\mathrm{R}p_{2}}_{*}(\mathscr{E}^{\mathrm{R}\hbar}
  \stackrel{\mathrm{L}}\otimes_{p_{1}^{-1}\mathscr{O}_{X}^{\hbar}}p_{1}^{-1}
  (\bigoplus\mathscr{M}_{i})^{cc})^{cc}.$
  \item [] $\simeq{\mathrm{R}p_{2}}_{*}\mathrm{R}\mathscr{H}om_{p_{2}^{-1}\mathscr{O}_{Y}^{\hbar}}
  (p_{2}^{-1}\mathscr{O}_{Y}^{\mathrm{loc}}/\mathscr{O}_{Y}^{\hbar}[-1],
  \mathscr{E}^{\mathrm{R}\hbar}\stackrel{\mathrm{L}}\otimes_{p_{1}^{-1}\mathscr{O}_{X}^{\hbar}}
  p_{1}^{-1}\mathrm{R}\mathscr{H}om_{\mathscr{O}_{X}^{\hbar}}
  (\mathscr{O}_{X}^{\mathrm{loc}}/\mathscr{O}_{X}^{\hbar}[-1],\bigoplus\mathscr{M}_{i}))$
\end{itemize}
and
\begin{itemize}
  \item [] ${\mathrm{R}p_{2}}_{*}(\bigoplus(\mathscr{E}^{\mathrm{R}\hbar}
  \stackrel{\mathrm{L}}\otimes_{p_{1}^{-1}\mathscr{O}_{X}^{\hbar}}p_{1}^{-1}
  \mathscr{M}_{i}))^{cc}$
  \item [] $\simeq{\mathrm{R}p_{2}}_{*}\mathrm{R}\mathscr{H}om_{p_{2}^{-1}\mathscr{O}_{Y}^{\hbar}}
  (p_{2}^{-1}\mathscr{O}_{Y}^{\mathrm{loc}}/\mathscr{O}_{Y}^{\hbar}[-1],
  \mathscr{E}^{\mathrm{R}\hbar}\stackrel{\mathrm{L}}\otimes_{p_{1}^{-1}\mathscr{O}_{X}^{\hbar}}
  p_{1}^{-1}\bigoplus\mathscr{M}_{i}))$.
\end{itemize}
From the natural morphism
\begin{center}
$\bigoplus(\mathscr{M}_{i})^{cc}\rightarrow(\bigoplus\mathscr{M}_{i})^{cc}$ and
$\mathscr{M}_{i}\simeq(\mathscr{M}_{i})^{cc}$,
\end{center}
we get the following morphism
\begin{center}
${\mathrm{R}p_{2}}_{*}(\bigoplus(\mathscr{E}^{\mathrm{R}\hbar}
  \stackrel{\mathrm{L}}\otimes_{p_{1}^{-1}\mathscr{O}_{X}^{\hbar}}p_{1}^{-1}
  \mathscr{M}_{i}))^{cc}\rightarrow{\mathrm{R}p_{2}}_{*}(\mathscr{E}^{\mathrm{R}\hbar}
  \stackrel{\mathrm{L}}\otimes_{p_{1}^{-1}\mathscr{O}_{X}^{\hbar}}p_{1}^{-1}
  (\bigoplus\mathscr{M}_{i})^{cc})^{cc}.$
\end{center}
Taking $\mathrm{gr}_{\hbar}$ functor and use Corollary 6.15 and Proposition 6.19, we get
\begin{center}
$\overline{\bigoplus}(\Phi_{\mathscr{E}^{\mathrm{R}\hbar}}\mathscr{M}_{i})
\simeq\Phi_{\mathscr{E}^{\mathrm{R}\hbar}}(\overline{\bigoplus}
\mathscr{M}_{i})$.
\end{center}
$\hfill \square$\\

\noindent {\bf{Proposition 6.22.}} \emph{The functor $\Phi_{\mathscr{E}^{\mathrm{R}\hbar}}
(\mathscr{M}):={\mathrm{R}p_{2}}_{*}(\mathscr{E}^{\mathrm{R}\hbar}
\stackrel{\mathrm{L}}\otimes_{p_{1}^{-1}\mathscr{O}_{X}^{\hbar}}p_{1}^{-1}
\mathscr{M})^{{cc}}$ is fully faithful sending $\mathrm{D}^{\mathrm{b}}_{\mathrm{coh}}
(\mathscr{O}_{X}^{\hbar})$ to $\mathrm{D}^{\mathrm{b}}_{\mathrm{coh}}
(\mathscr{O}_{Y}^{\hbar})$ and has a right adjoint functor.}\\

\noindent
{\emph{{Proof}}}. Consider the following morphism
\begin{center}
$\mathrm{RHom}(\mathscr{M},\mathscr{N})\longrightarrow
\mathrm{RHom}(\Phi_{\mathscr{E}^{\mathrm{R}\hbar}}\mathscr{M},
\Phi_{\mathscr{E}^{\mathrm{R}\hbar}}\mathscr{N})$.
\end{center}
Taking $\mathrm{gr}_{\hbar}$ functor, and using $\Phi_{\mathscr{E}}$ is fully faithful, we get
$\Phi_{\mathscr{E}^{\mathrm{R}\hbar}}$ is fully faithful.
Now use Theorem 6.5, Proposition 6.13 and Proposition 6.21 and the finiteness of
convolution.$\hfill \square$\\

\noindent {\bf{Proposition 6.23.}} ([Huy06, Lemma 1.50]). \emph{Let $F:\mathcal{D}\longrightarrow
\mathcal{D}'$ be a fully faithful exact functor between triangulated
categories and suppose that $F$ has a right adjoint $F\dashv{H}$.
Then $F$ is an equivalence if and only if for any $C\in\mathcal{D}'$
the triviality of $H(C)$, i.e., $H(C)\simeq{0}$ implies $C\simeq{0}$.}\\

\noindent {\bf{Theorem 6.24.}} \emph{The following functor
\begin{center}
$\Phi_{\mathscr{E}^{\mathrm{R}\hbar}}:\mathrm{D}_{\mathrm{qcc}}(\mathscr{O}_{X}^{\hbar})
\rightarrow\mathrm{D}_{\mathrm{qcc}}(\mathscr{O}_{Y}^{\hbar})$
\end{center}
is an equivalence.}\\

\noindent
{\emph{{Proof}}}. Fully faithfulness of $\Phi_{\mathscr{E}^{\mathrm{R}\hbar}}$ is proved
in Proposition 6.22.

Now we prove the surjectivity. Denote by $G$ the right adjoint functor of
$\Phi_{\mathscr{E}^{\mathrm{R}\hbar}}$.
Then for $\mathscr{M}\in\mathrm{D}_{\mathrm{qcc}}(\mathscr{O}_{X}^{\hbar})$ and
$\mathscr{N}\in\mathrm{D}_{\mathrm{qcc}}(\mathscr{O}_{Y}^{\hbar})$, we have
\begin{center}
$\mathrm{RHom}(\Phi_{\mathscr{E}^{\mathrm{R}\hbar}}
\mathscr{M},\mathscr{N})\simeq\mathrm{RHom}(\mathscr{M},G\mathscr{N})$.
\end{center}
If $G\mathscr{N}\simeq{0}$, then $\mathrm{RHom}(\mathrm{gr}_{\hbar}
\Phi_{\mathscr{E}^{\mathrm{R}\hbar}}\mathscr{M},\mathrm{gr}_{\hbar}
\mathscr{N})\simeq\mathrm{RHom}(\Phi_{\mathscr{E}}\mathrm{gr}_{\hbar}
\mathscr{M},\mathrm{gr}_{\hbar}\mathscr{N})\simeq{0}$.
Denote by $F'$ the inverse functor of $F$. Let
$\mathscr{M}=(F'\mathscr{M}')^{\mathrm{R}\hbar}$
with $\mathscr{M}'\in\mathrm{D}^{\mathrm{b}}_{\mathrm{coh}}
(\mathscr{O}_{Y})$ which implies that $\mathrm{gr}_{\hbar}\mathscr{N}\simeq{0}$
by Example 6.6. Then $\mathscr{N}\simeq{0}$ since $\mathscr{N}$ is cohomologically
complete.

This proves the equivalence
\begin{center}
$\mathrm{D}_{\mathrm{qcc}}(\mathscr{O}_{X}^{\hbar})\stackrel{\sim}
\longrightarrow\mathrm{D}_{\mathrm{qcc}}(\mathscr{O}_{Y}^{\hbar})$
\end{center}
by Proposition 6.23.$\hfill \square$\\

\noindent {\bf{Theorem 6.25.}} \emph{The equivalent functor
$\Phi_{\mathscr{E}^{\mathrm{R}\hbar}}$ induces the following
equivalence (we keep the same notation)}
\begin{center}
$\Phi_{\mathscr{E}^{\mathrm{R}\hbar}}:\mathrm{D}^{\mathrm{b}}_{\mathrm{coh}}
(\mathscr{O}_{X}^{\hbar})\stackrel{\sim}\longrightarrow
\mathrm{D}^{\mathrm{b}}_{\mathrm{coh}}(\mathscr{O}_{Y}^{\hbar})$.
\end{center}
\emph{{Proof}}. By Theorem 6.24, we only need to prove the surjectivity.
Let $\mathscr{N}\in\mathrm{D}^{\mathrm{b}}_{\mathrm{coh}}(\mathscr{O}_{Y})$.
Then there exists $\mathscr{M}\in\mathrm{D}_{\mathrm{qcc}}(\mathscr{O}_{X}^{\hbar})$
such that $\Phi_{\mathscr{E}^{\mathrm{R}\hbar}}(\mathscr{M})\simeq\mathscr{N}$.
Taking $\mathrm{gr}_{\hbar}$ functor, we get $\Phi_{\mathscr{E}}(\mathrm{gr}_{\hbar}(\mathscr{M}))
\simeq\mathrm{gr}_{\hbar}\mathscr{N}\simeq{F}(\mathscr{F})\simeq\Phi_{\mathscr{E}}
(\mathscr{F})$ by the equivalence of $F$ for some $\mathscr{F}\in
\mathrm{D}^{\mathrm{b}}_{\mathrm{coh}}(\mathscr{O}_{X})$. Since
$\Phi_{\mathscr{E}}$ is an equivalence, this implies that $\mathrm{gr}_{\hbar}
\mathscr{M}\simeq\mathscr{F}\in\mathrm{D}^{\mathrm{b}}_{\mathrm{coh}}
(\mathscr{O}_{X})$ and hence $\mathscr{M}\in\mathrm{D}^{+}_{\mathrm{qcc}}
(\mathscr{O}_{X}^{\hbar})$ by Proposition 6.12. Since
$\mathrm{gr}_{\hbar}\mathscr{M}\in\mathrm{D}^{\mathrm{b}}_{\mathrm{coh}}
(\mathscr{O}_{X}),\mathscr{M}\in\mathrm{D}^{+}_{\mathrm{coh}}
(\mathscr{O}_{X}^{\hbar})$ by Proposition 6.13. Hence by
Lemma 2.1 and Nakayama's Lemma, we obtain $\mathscr{M}\in\mathrm{D}^{\mathrm{b}}_{\mathrm{coh}}
(\mathscr{O}_{X}^{\hbar})$. Thus, $\Phi_{\mathscr{E}^{\mathrm{R}\hbar}}
(\mathscr{M})\simeq\mathscr{N}$ as desired.$\hfill \square$

\begin{center}

\end{center}

\noindent
{Institute of Mathematics, Academia Sinica, Taipei 106, Taiwan}\\
E-mail address: houyi@math.sinica.edu.tw; d93221002@ntu.edu.tw

\end{document}